\documentclass{article}
\usepackage{graphicx,amsmath,amsthm,amssymb,float} % Required for inserting images

\newtheorem{theorem}{Theorem}[section]
\newtheorem{definition}[theorem]{Definition}
\newtheorem{remark}[theorem]{Remark}

\newtheorem{question}[theorem]{Question}

\usepackage{xcolor}
\usepackage{pgf,tikz}

\usepackage{hyperref}
\usepackage{lmodern}
\usepackage{microtype}
\usepackage{amsmath,amssymb,mathtools}
\usepackage{enumitem}
\usepackage{graphicx}
\usepackage{float}
\usepackage[h]{esvect}

\usepackage{chngcntr}

\makeatletter
\newcommand{\customitem}[1]{%
\item[#1]\protected@edef\@currentlabel{#1}%
}
\makeatother

\usepackage{hhline}
\usepackage{diagbox}
\usepackage{multirow}

\definecolor{myred}{rgb}{0.81, 0.06, 0.13}
\definecolor{myblue}{rgb}{0,0.5,0.6}

\def\x{\mathbf{x}}
\def\y{\mathbf{y}}

\def\p{\mathbf{p}}

\def\R{{\mathbb R}}
\def\Z{{\mathbb Z}}

\theoremstyle{plain}

\theoremstyle{definition}

\newcommand{\HH}{{\mathcal H}}
\renewcommand{\H}{\HH^1}

\usepackage[margin=2cm]{geometry}

\usetikzlibrary{arrows.meta,calc,decorations.pathreplacing,patterns}

\setlist[itemize]{leftmargin=1.6em,itemsep=0.45em,topsep=0.4em}
\setlist[enumerate]{leftmargin=1.8em,itemsep=0.55em,topsep=0.4em}
\newcommand{\norm}[1]{\left\lVert #1\right\rVert}
\newcommand{\potential}[2]{a_{#1,#2}}

\begin{document}
\title{Open problems from the Length Optimization event at MDS 2026 \\6-7 July 2026}
\author{Yana Teplitskaya \footnote{You can write to Yana to get in touch with any coauthor: \href{teplitskayayana@gmail.com}{teplitskayayana@gmail.com}}}
\date{Institute of Mathematics and Informatics Bulgarian Academy of Sciences, Sofia, Bulgaria}
\maketitle
%\section{Main}

During the Minisymposium on Length Optimization as part of the International Conference “Mathematics Days in Sofia – 2026”, the participants collected the following open problems:
 
\begin{itemize}
    \item Find any non-trivial upper and lower bounds for the degree of a branching point in the Gilbert-Steiner problem with cost-function $c(x) = x^p$ in the metric space $X = \mathbb{R}^d$ for large $d$ and fixed $p > 1/2$ (see Section~\ref{DC}).
    \item Determine whether a typical trajectory of a given branched optimal traffic plan has bounded length (see Section~\ref{AP}). To find the optimal transport plan between the Dirac measure in the center of a circle and the Lebesgue measure on the circle (see Section~\ref{AP}).
    \item To find the maximal distance minimizers for any explicit compact sets (see Section~\ref{YT}).
    \item What is the minimal $C$ such that for any set of $n$ points in the three-dimensional unit cube, there is a tour through all these $n$ points with the sum of the cubes of its edge lengths at most $C$ (see Section~\ref{LG}).
    \item To explore parallels between caustics and Steiner trees: do caustics really need convexity (see Section~\ref{MS}).
\end{itemize}

\section{Question of Danila Cherkashin \href{mailto:jiocb@math.bas.bg}{jiocb@math.bas.bg}}
\label{DC}

Let us define the following Gilbert--Steiner problem.

\begin{definition}
    Let $\mu^+,\mu^-$ be two finite measures on a metric space $(X,\rho(\cdot,\cdot))$ with finite supports such that total masses $\mu^+(X)=\mu^-(X)$ are equal. Let $V\subset X$  be a finite set containing the support of the signed measure $\mu^+-\mu^-$; the elements of $V$ are called \emph{vertices}. Further, let $E$ be a finite collection of unordered pairs $\{x,y\}\subset V$ which we call \emph{edges}. So, $(V,E)$ is a  simple undirected finite graph. Assume that for every $\{x,y\}\in E$ two non-zero real numbers $m(x,y)$ and $m(y,x)$ are defined so that $m(x,y)+m(y,x)=0$. This data set is called a \emph{$(\mu^+,\mu^-)$-flow} if 
    \[
    \mu^+ - \mu^- = \sum_{\{x,y\}\in E} m(x,y)\cdot (\delta_y-\delta_x)
    \]
    where $\delta_x$ denotes a delta-measure at $x$ (note that the summand $m(x,y)\cdot (\delta_y-\delta_x)$ is well-defined in the sense that it does not depend on the order of $x$ and $y$).
\end{definition}

Let $c \colon [0,\infty)\to [0,\infty)$ be a \emph{cost function}. The expression
\[
\sum_{\{x,y\}\in E} c(|m(x,y)|) \cdot \rho(x,y)
\]
is called the \emph{Gilbert functional} of the $(\mu^+,\mu^-)$-flow.

The \textit{Gilbert--Steiner problem} is to find a flow which minimizes the Gilbert functional with cost function $c(x) = x^p$, for a fixed $p \in (0,1)$, in the metric space $X = \mathbb{R}^d$; we call such a solution a \emph{minimal flow}.

\begin{question}
    Find any non-trivial upper and lower bounds for the maximum possible degree of a branching point for large $d$ and fixed $p > 1/2$.
\end{question}

One can easily reduce this problem to a neighborhood of a branching point (and by the hereditary property the restricted flow will also be minimal). Then, for a small enough neighborhood, there is a structural theorem in the case of a single source (see also~\cite{cherkashin2025branching}). Actually, this theorem works for a monotone, concave and non-negative cost function $c(x)$ such that $c(0) = 0$.

\begin{theorem}[Volz--Brazil--Ras--Swanepoel--Thomas, 2013~\cite{volz2013gilbert}] \label{th:volz}
Consider a single source at $q$ and sinks $p_1,\dots,p_n \in \mathbb{R}^d$. 
Let the mass demand at $p_i$ be $m_i > 0$ for $i=1,\dots,n$. 
Then a star flow with a unique branching point at the origin $o$ is optimal if and only if
\begin{equation}\label{eq:wstp}
\sum_{i=1}^n c(m_i) e_i + c\left( \sum_{i=1}^n m_i \right) e = \bar{0}
\end{equation}
and
\begin{equation}\label{eq:split1}
\left \langle \sum_{i\in I} c(m_i)e_i,\sum_{i\in I} c(m_i)e_i \right \rangle \leq c^2 \left (\sum_{i\in I} m_i \right ) \quad \text{ for all } I\subseteq\{1,\dots,n\},
\end{equation}
where $e_i$ is the unit vector codirected with $op_i$ and $e$ is the unit vector codirected with $oq$.
\end{theorem}

To obtain an upper bound for a fixed $p > 1/2$ note that the angle between every pair of edges is bounded from below by some $\alpha(p)$. 
Thus one can use a bound on the cardinality of a spherical code $A(d,\alpha(p))$, see~\cite{BDM12} for the survey and~\cite{kabatiansky1978bounds} for the best asymptotic bound. However numerical simulations show that for certain configurations a critical set $I$ may have cardinality greater than the dimension $d$.

In the plane the branching is always triple and this does not depend on $p$. A major question is whether there is a bound that depends only on $d$ for every $d > 2$?

\bibliographystyle{abbrv}
\bibliography{main}

\section{Questions of Antoine Prouff}
\label{AP}

\subsection{Uniform length bound for optimal irrigation plans}

Here we discuss Problem 15.13 from~\cite{BCM}. Let $\alpha \in (0, 1)$ be fixed\footnote{this parameter is actually $p$ from the previous discrete setting} and $\mu^-, \mu^+$ be two mutually singular and compactly supported probability measures on $\mathbb{R}^d$. Let $\mathbf{P}$ be an optimal traffic plan connecting $\mu^-$ to $\mu^+$, namely a probability measure $\mathbf{P}$ on the space $\text{Lip}$ of {$1$-Lipschitz} curves $\gamma : [0, \infty) \to \mathbb{R}^d$, with minimal cost
\begin{equation*}
\mathcal{E}^\alpha(\mathbf{P})
	= \int_{\text{Lip}} \left( \int_0^\infty \theta_{\mathbf{P}}^{\alpha-1}\left(\gamma(t)\right) |\dot{\gamma}(t)| \, dt \right) \, d\mathbf{P}(\gamma) ,
\end{equation*}
where $\theta_{\mathbf{P}}(x) = \mathbf{P}\left(\{\gamma \in \text{Lip} \, : \, x \in \gamma([0, \infty)) \} \right)$.
Assume that $\mathcal{E}^\alpha(\mathbf{P}) < \infty$ (this is automatic if $\alpha > 1 - 1/d$). Since $\theta_{\mathbf{P}} \le 1$, we have
\begin{equation*}
\int_{\text{Lip}} \text{length}(\gamma) \, d\mathbf{P}(\gamma)
	= \int_{\text{Lip}} \left( \int_0^\infty |\dot{\gamma}(t)| \, dt \right) \, d\mathbf{P}(\gamma)
	\le \mathcal{E}^\alpha(\mathbf{P})
	< \infty ,
\end{equation*}
namely the average length of trajectories is finite, or in other words, the length function belongs to $L^1(\text{Lip}, \mathbf{P})$. In particular, Markov's inequality yields
\begin{equation*}
\mathbf{P}\left(\{\text{length}(\gamma) \ge L\}\right)
	\le \dfrac{c}{L} ,
		\qquad \forall L \ge 1 .
\end{equation*}
The problem that we are asking is whether there exists a constant $c_\infty > 0$ such that
\begin{equation} \label{e:Linfty}
\text{length}(\gamma) \le c_\infty,
	\qquad \text{for {$\mathbf{P}$-almost} every $\gamma$,}
\end{equation}
or equivalently, whether the length function belongs to $L^\infty(\text{Lip}, \mathbf{P})$.
An intermediate goal would be to investigate the possibility of the length being in $L^p(\text{Lip}, \mathbf{P})$, or the slightly weaker estimate
\begin{equation*}
\mathbf{P}\left(\{\text{length}(\gamma) \ge L\}\right)
	\le \dfrac{C}{L^p} ,
		\qquad \forall L \ge 1 .
\end{equation*}

Reasoning by contradiction, the existence of a positive measure set of arbitrarily long trajectories implies the presence of ``quasi-cycles" in the sense of~\cite{CDRMPP}, namely curves that almost close up after traveling a certain distance (recall that {$\mathbf{P}$-almost} every curve lives in the convex hull of $\text{supp} \, \mu^- \cup \text{supp} \, \mu^+$ by optimality and finiteness of the cost). Unfortunately, for the problem~\eqref{e:Linfty}, the measure/weight of those quasi-cycles does not seem to be enough to make a short-cutting argument work (i.e.\ construct a competitor by removing the quasi-cycle and reconnect locally to keep the marginal condition satisfied, see~\cite{CDRM,CDRMPP}). Indeed, it is not obvious that the energy saved by short-cutting the quasi-cycle is enough to balance the additional cost of the ``local reconnection". Perhaps relaxing the problem to an $L^p$ bound would give us a little more room to build a competitor using this short-cutting procedure.

One particular case to keep in mind is when the supports of $\mu^-$ and $\mu^+$ are disjoint and $\mu^-$ is a Dirac delta. Then the multiplicity $\theta_{\mathbf{P}}$ is non-increasing along almost every trajectory. Understanding the decay rate of the multiplicity could prove useful to make progress towards~\eqref{e:Linfty}.

\subsection{Irrigating the unit circle from the origin}

We wish to describe the set of optimal traffic plans in a very simple and specific situation.
In $\mathbb{R}^2$, set $\mu^- = \delta_0$ and $\mu^+ = \frac{1}{2\pi} \mathcal{H}^1 \, {\mathbin{\vrule height 1.6ex depth 0pt width 0.13ex\vrule height 0.13ex depth 0pt width 1.3ex}} \, \mathbb{S}^1$. Let $\mathbf{P} = \mathbf{P}_\alpha$ be an optimal traffic plan for the {$\alpha$-energy}, $\alpha \in (0, 1)$, connecting $\mu^-$ to $\mu^+$.
\begin{itemize}
\item Is there branching at the origin? If yes, what is the degree as a function of $\alpha$? The local picture seems to be related to Danila's question.
\item In connection with numerical simulations: Does $\mathbf{P}$ always have symmetries, or can we come up with an instance of an optimal traffic plan without obvious symmetries?
\item What about the number of branching points in $B_r(0)$ as $r \to 1^-$? Can we quantify the decay of the multiplicity $\theta_{\mathbf{P}}$ along trajectories?
\end{itemize}

Notice that optimizers are known explicitly for the end-point cases $\alpha \in \{0, 1\}$. As a starting point, we could investigate whether we can obtain any interesting results from a continuity argument as $\alpha \to 1^-$.

\newcommand{\etalchar}[1]{$^{#1}$}

\section{Question of Yana Teplitskaya}
\label{YT}
For a given compact set $M \subset \mathbb{R}^d$ and a given positive number $r > 0$, we seek a connected compact set $\Sigma$ of minimal length (one-dimensional Hausdorff measure $\H$) such that $M \subset \overline{B_r(\Sigma)}$ where $\overline{B_r(\Sigma)}$ denotes the closed $r$-neighborhood of the desired set $\Sigma$.  
We will call such a solution $\Sigma$ a
\textit{Maximal distance minimizer} (\textit{MDM} or $r$-minimizer).

A solution always exists and contains no loops (homeomorphic image of a circle), see~\cite{paolini2004qualitative}.
For detailed information on MDMs we refer the reader to our review~\cite{cherkashin2022overview}.

\begin{question}
    What is a maximal distance minimizer for a circle of radius $R$ ($M=\partial B_R$) with $R>r$? This remains open for $4.98r>R>r$.
\end{question}

\begin{figure}
  \begin{center}
  %\begin{figure}[h]
%	\begin{center}
		\definecolor{yqqqyq}{rgb}{0.,0.,1.}%{0.5019607843137255,0.,0.5019607843137255}
		\definecolor{xdxdff}{rgb}{0.,0.,0.}
		\definecolor{ffqqqq}{rgb}{0.,0.,1.}
    \definecolor{ttzzqq}{rgb}{0.,0.,0.}
    \begin{tikzpicture}[line cap=round,line join=round,x=1.0cm,y=1.0cm,scale=1.5]
    \clip(-5.1,-2.15198488708635626) rectangle (0.1,2.2268250892462824);
    \draw [very thick, color=ttzzqq] (-3.,2.)-- (-2.,2.);
    \draw [very thick, color=ttzzqq] (-2.,-2.)-- (-3.,-2.);
    \draw [very thick, color=ffqqqq] (-3.,-1.)-- (-2.,-1.);
    \draw [dotted] (-3.,1.)-- (-2.,1.);
    \draw [shift={(-3.,0.)},very thick, color=ttzzqq]  plot[domain=1.5707963267948966:4.71238898038469,variable=\t]({1.*2.*cos(\t r)+0.*2.*sin(\t r)},{0.*2.*cos(\t r)+1.*2.*sin(\t r)});
    \draw [shift={(-3.,0.)},dotted]  plot[domain=1.5707963267948966:4.71238898038469,variable=\t]({1.*1.*cos(\t r)+0.*1.*sin(\t r)},{0.*1.*cos(\t r)+1.*1.*sin(\t r)});
    \draw [shift={(-2.,0.)},dotted]  plot[domain=-1.5707963267948966:1.5707963267948966,variable=\t]({1.*1.*cos(\t r)+0.*1.*sin(\t r)},{0.*1.*cos(\t r)+1.*1.*sin(\t r)});
    \draw [shift={(-2.,0.)},very thick,color=ttzzqq]  plot[domain=-1.5707963267948966:1.5707963267948966,variable=\t]({1.*2.*cos(\t r)+0.*2.*sin(\t r)},{0.*2.*cos(\t r)+1.*2.*sin(\t r)});
    \draw [dash pattern=on 2pt off 2pt] (-2.5008100281931047,2.) circle (1.cm);
    \draw [very thick,color=yqqqyq] (-3.7309420990521054,0.6824394829091469)-- (-3.1832495111022516,1.2690579009478955);
    \draw [very thick,color=yqqqyq] (-1.8178907146013357,1.2695061868000874)-- (-1.2695061868000874,0.6829193135917684);
    \draw [shift={(-2.,0.)},very thick,color=ffqqqq]  plot[domain=-1.5707963267948966:0.7517515639553677,variable=\t]({1.*1.*cos(\t r)+0.*1.*sin(\t r)},{0.*1.*cos(\t r)+1.*1.*sin(\t r)});
    \draw [shift={(-3.,0.)},very thick,color=ffqqqq]  plot[domain=2.390497746093128:4.771687465439039,variable=\t]({1.*1.*cos(\t r)+0.*1.*sin(\t r)},{0.*1.*cos(\t r)+1.*1.*sin(\t r)});
    \begin{scriptsize}
    \draw[color=ttzzqq] (-4.8089822121844977,0.2025724503290166) node {$M$};
    %\draw [fill=xdxdff] (-2.5008100281931047,2.) circle (1.5pt);
    %\draw [fill=ffqqqq] (-3.1832495111022516,1.2690579009478955) circle (1.5pt);
    %\draw [fill=ffqqqq] (-1.8178907146013357,1.2695061868000874) circle (1.5pt);
    \draw[color=blue] (-3.763965847825129,-0.1818662773850787) node {$\Sigma$};
   % \draw[color=black] (-2.059620164509451,0.8028001329511188) node {$M_r$};
    \end{scriptsize}
    \end{tikzpicture}
   
	%\end{center}
  %  \end{figure}
  \caption{A horseshoe for a closed curve with a big radius of curvature}
  \label{fig:horseshoe}
  \end{center}
\end{figure}
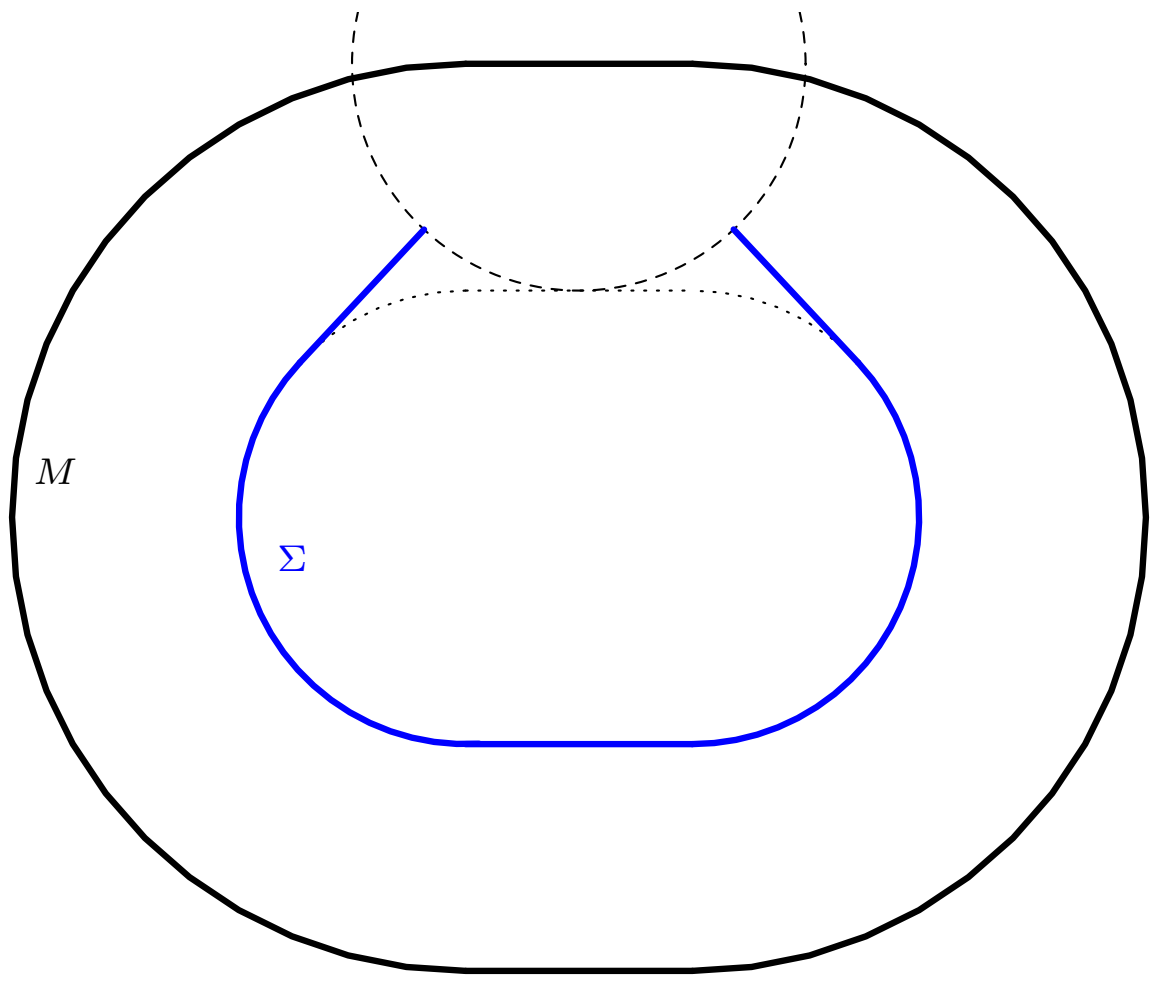

Miranda, Paolini and Stepanov conjectured in~\cite{miranda2006one} that any minimizer for a circle of radius $R > r$ is a horseshoe, i.e. a union of an arc of a circle of radius $R-r$ and two segments of its tangent lines.

The conjecture is confirmed in~\cite{cherkashin2018horseshoe} for $R>4.98r$. For smaller values with $4.98r>R > r$ the hypothesis remains unproved.

\begin{question}
     What is the set of maximal distance minimizers when the input set $M$ is an $R$-stadium (which is the boundary of the $R$-neighborhood of a segment)? This is open for $r<R<5r$.
\end{question}

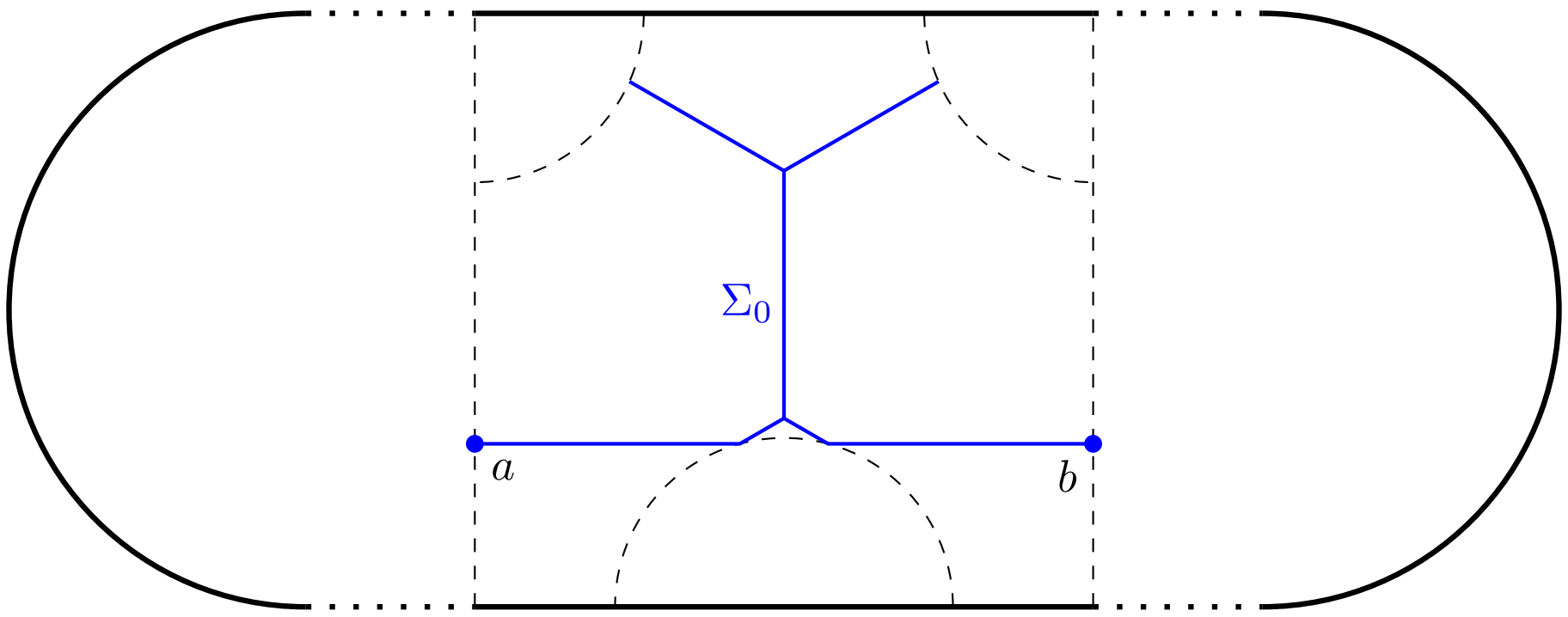
\begin{figure}
  \begin{center}
  %\begin{figure}[h]
 %   \centering

\begin{tikzpicture}[scale=1.3]

\foreach \y in {0, 3.51389977019889} {
    \draw[very thick, loosely dotted] 
        (-1, \y) -- (0, \y);
    \draw[very thick]
        (0, \y) -- (3.66, \y);
    \draw[very thick, loosely dotted] 
        (3.66, \y) -- (4.66, \y);
}
\draw[very thick] (-1, 0) arc (270:90:1.757);
\draw[very thick] (4.66, 0) arc (-90:90:1.757);

%\draw[thin, dashed] (0.2614137092297641 + 1, 0.0) arc (0:180:1);
%\draw[thin, dashed] (0.7842411276892922 + 1, 0.0) arc (0:180:1);
%\draw[thin, dashed] (1.3070685461488205 + 1, 0.0) arc (0:180:1);
\draw[thin, dashed] (1.8298959646083486 + 1, 0.0) arc (0:180:1);
%\draw[thin, dashed] (2.3527233830678767 + 1, 0.0) arc (0:180:1);
%\draw[thin, dashed] (2.875550801527405 + 1, 0.0) arc (0:180:1);
%\draw[thin, dashed] (3.398378219986933 + 1, 0.0) arc (0:180:1);
%\draw[thin, dashed] (1.8298959646083486 + 1, 3.75) arc (0:-180:1);

\draw[blue, thick] 
(0.0, 0.9652268503449002) -- 
(1.5684822553785844, 0.9652268503449002) --
(1.8298959646083486, 1.11615412573856) --
(2.0913096738381127, 0.9652268503449002) -- 
(3.659791929216697, 0.9652268503449002);

\node[blue] at (1.61,1.8) {$\Sigma_0$};

\draw[dashed] (0.0, 0.0) -- (0.0, 3.51389977019889);
\draw[dashed] (3.659791929216697, 0.0) -- (3.659791929216697, 3.51389977019889);

\fill[blue] (0.0, 0.9652268503449002) 
    circle (1.5pt) node[black, below right] {$a$};
\fill[blue] (3.659791929216697, 0.9652268503449002) 
    circle (1.5pt) node[black, below left] {$b$};

\draw[blue, thick] 
(1.8298959646083486, 1.11615412573856) -- 
(1.8298959646083486, 2.58208251632886);

\draw[thin, dashed] (1, 3.51389977019889) arc (0:-90:1);
\draw[thin, dashed] (2.66, 3.51389977019889) arc (180:270:1);

\draw[blue, thick]
(1.8298959646083486, 2.58208251632886) -- 
(0.9149479823041743, 3.11032798020668);

\draw[blue, thick] 
(1.8298959646083486, 2.58208251632886) -- 
(2.744843946912523, 3.11032798020668);

\end{tikzpicture}

   % \caption{Horseshoe is not a minimizer for long enough stadium with $R < 1.75r$.}
  %  \label{stadion}
%\end{figure}
  \caption{A competitor shorter than the horseshoe for a stadium} 
  \label{fig:stadium}
  \end{center}
\end{figure}

 In~\cite{cherkashin2018horseshoe} it is proved that for a closed curve with minimal radius of curvature $R>5r$ any MDM is a horseshoe, i.e. a union of an arc of a curve parallel to an arc of $M$ (and hence having a radius of curvature at least $R-r$) and two segments of its tangent lines (see Figure~\ref{fig:horseshoe}). Note that if a curve has a strictly positive radius of curvature then it is at least $C^{1,1}$-smooth; for such curves the parallel (offset) curve is well-defined by the map $s \mapsto \gamma(s) + r\vec{n}_i(s)$, where $\vec{n}_i(s)$ is the inner unit normal vector to a given curve $M=\gamma(s)$. 
Thus, for $R>5r$, every solution for the $R$-stadium is a horseshoe. It is also known that, if $r<R<1.75r$ and the segment defining
the $R$-stadium is sufficiently long, then no MDM is a horseshoe (see~\cite{cherkashin2022overview} and Figure~\ref{fig:stadium}).
%It is also known (see~\cite{cherkashin2022overview} and Figure~\ref{fig:stadium}) that for $R<1.75r$ the solution for an $R$-stadium is not a horseshoe. 

 Thus the following subquestions remain open:
\begin{itemize}
    \item To find the supremum of $C$ such that there exists an MDM for a $Cr$-stadium which is not a horseshoe.
    \item To find the infimum of $c$ such that for any closed curve with minimal radius of curvature $cr$ every MDM is a horseshoe.
\end{itemize}
It might also be helpful to know that any minimizer for an $R$-stadium consists of line segments and of curves parallel to the curves contained in an $R$-stadium (these curves coincide with arcs of circles of radius $R-r$).

\begin{question}
  What is the set of MDMs for any explicit set you like?
 \end{question}

  For example for 
    \begin{itemize}
        \item A disc or ball $B_R$ in $\mathbb{R}^d$, $d\geq 2$. 
        A very challenging question, nothing is known except degenerate cases and the usual regularity properties such as absence of loops.
        \item A rectangle (the boundary). For a small enough $r>0$ the question is answered (and is surprisingly asymmetric), see~\cite{cherkashin2021maximal}.
    \end{itemize}

\section{Question of Alexey Gordeev}

\label{LG}
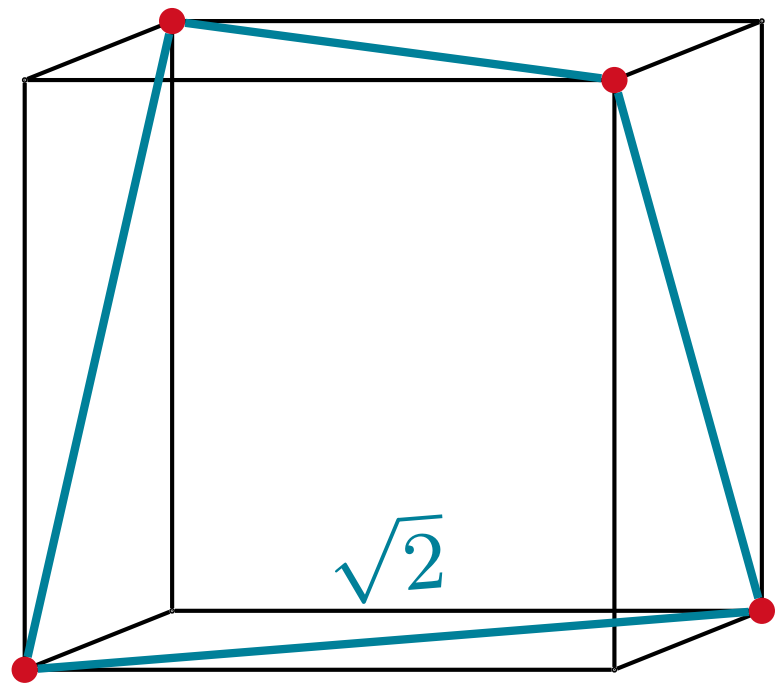
\begin{figure}[!ht]
\begin{center}
\resizebox{0.3\textwidth}{!}{%
\tikzset{cubevert/.style={circle,draw,fill,inner sep=0,minimum size=0.005cm}, vertex/.style={circle,draw=myred,fill=myred,inner sep=0pt,minimum size=0.075cm}, every edge/.style={draw}}
\begin{tikzpicture}
  
    \node[cubevert] (000) at (0,0) {};
    \node[cubevert] (010) at (2,0) {};
    \node[cubevert] (001) at (0,2) {};
    \node[cubevert] (011) at (2,2) {};
    \node[cubevert] (100) at (0.5,0.2) {};
    \node[cubevert] (110) at (2.5,0.2) {};
    \node[cubevert] (101) at (0.5,2.2) {};
    \node[cubevert] (111) at (2.5,2.2) {};
    
    \draw (000) edge (001)
          (010) edge (011)
          (100) edge (101)
          (110) edge (111);
    \draw (000) edge (010)
          (001) edge (011)
          (100) edge (110)
          (101) edge (111); 
    \draw (000) edge (100)
          (001) edge (101)
          (010) edge (110)
          (011) edge (111);

    \node[vertex] (a) at (0,0) {};
    \node[vertex] (b) at (0.5,2.2) {};
    \node[vertex] (c) at (2,2) {};
    \node[vertex] (d) at (2.5,0.2) {};
  
    \draw (a) edge[myblue,thick] (b)
          (b) edge[myblue,thick] (c)
          (c) edge[myblue,thick] (d)
          (d) edge[myblue,thick] node[above,sloped]{\footnotesize $\sqrt{2}$} (a);
\end{tikzpicture}}
\end{center}
\caption{The set of 4 points with the shortest tour of weight $2^{7 / 2}$.}
\label{fig:ex}
\end{figure}

In 1992, Bollob{\'a}s and Meir~\cite[Theorem~3]{bollobasTravellingSalesmanProblem1992} showed that for every $k \geq 1$ there exists a constant $c_k$ such that, for any $n$ points in the $k$-dimensional unit cube $[0, 1]^k$, one can find a tour $x_1, \dots, x_n$ through these $n$ points with
\[
\sum_{i = 1}^n |x_i - x_{i + 1}|^k \leq c_k,
\]
where $x_{n + 1} = x_1$ and $|x - y|$ is the Euclidean distance between $x$ and $y$.
Remarkably, this bound does not depend on $n$, the number of points.
They conjectured that the optimal constant is $c_k = 2 \cdot k^{k / 2}$ and showed that it cannot be taken lower than that.
This conjecture was recently revised for $k = 3$ by Balogh, Clemen and Dumitrescu~\cite[Theorem~1.3(i)]{baloghTravelingSalesmanProblem2024}, who showed that $c_3 \geq 4 \left( \sqrt{2} \right)^3 = 2^{7/2} > 2 \cdot 3^{3/2}$, see Figure~\ref{fig:ex}.
For $k = 2$, the conjecture is known to be true due to Newman~\cite[Problem 57]{newmanProblemSeminar1982}.
While it remains open for all $k > 2$, Gordeev~\cite[Theorem~1.7]{gordeevBollobasMeirTSPConjecture2026} recently showed that for large $k$ the conjecture holds asymptotically, i.e. $c_k = k^{k / 2} \cdot (2 + o_k(1))$.
The following question corresponds to, in a sense, the most aberrant and the least well-understood $3$-dimensional case of the conjecture.

\begin{question}
What is the minimal $C$ such that for any set of $n$ points $X \subseteq [0, 1]^3$ in the three-dimensional unit cube, there is a tour $x_1, \dots, x_n$ through these $n$ points with 
\[
\sum_{i = 1}^n |x_i - x_{i + 1}|^3 \leq C?
\]
Here, $x_{n + 1} = x_1$ and $|x - y|$ is the Euclidean distance between $x$ and $y$.
Is $C = 2^{7 / 2}?$    
\end{question}

% \bibliographystyle{abbrv}
% \bibliography{open_question_I}

\section{Question of Mikhail Shkolnikov}

\label{MS}
There are multiple parallels between Steiner trees and tropical caustics. One structural feature of the former is that convexity of the terminal locus is not required. Thus, one may expect that the current technical convexity assumption in tropical optics can be removed.

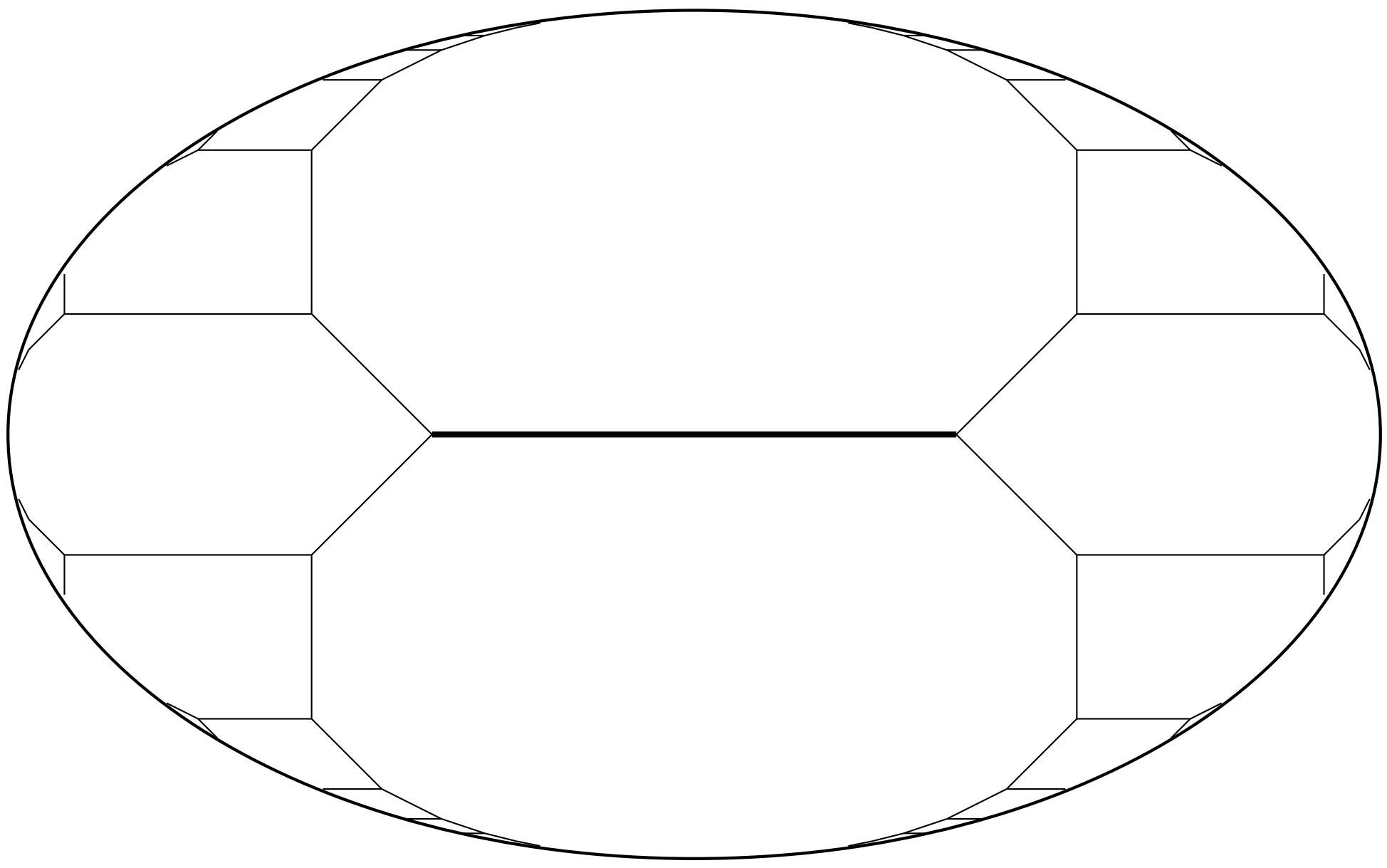
\begin{figure}[h!]\centering
\rotatebox{90}{
\begin{tikzpicture}
\begin{scope}[scale=3.94]
\begin{scope}[yscale= 1.6180339887498949]
\draw[thick] (0,0) circle (1);
\end{scope}
\draw[ultra thick](0.,0.61803398875)--(0.,-0.61803398875);
\draw(-0.959054704029,-0.41906747378)--(-0.970228164263,-0.36320017261);
\draw(-0.940424300357,-0.493589088465)--(-0.959054704029,-0.41906747378);
\draw(-0.940424300357,-0.493589088465)--(-0.940424300357,-0.547123155666);
\draw(-0.906338263824,-0.595847198066)--(-0.940424300357,-0.493589088465);
\draw(-0.906338263824,-0.595847198066)--(-0.906338263824,-0.677941159984);
\draw(-0.835969744015,-0.736584237683)--(-0.906338263824,-0.595847198066);
\draw(-0.835969744015,-0.736584237683)--(-0.835969744015,-0.875385672435);
\draw(-0.670440949108,-0.90211303259)--(-0.835969744015,-0.736584237683);
\draw(-0.670440949108,-1.16972939267)--(-0.715771340427,-1.12439900135);
\draw(-0.670440949108,-1.16972939267)--(-0.633478231341,-1.2436548282);
\draw(-0.670440949108,-0.90211303259)--(-0.670440949108,-1.16972939267);
\draw(-0.28407904384,-0.90211303259)--(-0.670440949108,-0.90211303259);
\draw(-0.28407904384,-1.48494113807)--(-0.3779862054,-1.48494113807);
\draw(-0.200037034439,-1.56898314747)--(-0.152208853468,-1.59289723796);
\draw(-0.28407904384,-1.48494113807)--(-0.200037034439,-1.56898314747);
\draw(-0.28407904384,-0.90211303259)--(-0.28407904384,-1.48494113807);
\draw(0.,-0.61803398875)--(-0.28407904384,-0.90211303259);
\draw(0.200037034439,-1.56898314747)--(0.152208853468,-1.59289723796);
\draw(0.28407904384,-1.48494113807)--(0.200037034439,-1.56898314747);
\draw(0.28407904384,-1.48494113807)--(0.3779862054,-1.48494113807);
\draw(0.28407904384,-0.90211303259)--(0.28407904384,-1.48494113807);
\draw(0.670440949108,-1.16972939267)--(0.633478231341,-1.2436548282);
\draw(0.670440949108,-1.16972939267)--(0.715771340427,-1.12439900135);
\draw(0.670440949108,-0.90211303259)--(0.670440949108,-1.16972939267);
\draw(0.835969744015,-0.736584237683)--(0.835969744015,-0.875385672435);
\draw(0.906338263824,-0.595847198066)--(0.906338263824,-0.677941159984);
\draw(0.940424300357,-0.493589088465)--(0.940424300357,-0.547123155666);
\draw(0.959054704029,-0.41906747378)--(0.970228164263,-0.36320017261);
\draw(0.940424300357,-0.493589088465)--(0.959054704029,-0.41906747378);
\draw(0.906338263824,-0.595847198066)--(0.940424300357,-0.493589088465);
\draw(0.835969744015,-0.736584237683)--(0.906338263824,-0.595847198066);
\draw(0.670440949108,-0.90211303259)--(0.835969744015,-0.736584237683);
\draw(0.28407904384,-0.90211303259)--(0.670440949108,-0.90211303259);
\draw(0.,-0.61803398875)--(0.28407904384,-0.90211303259);
\draw(0.959054704029,0.41906747378)--(0.970228164263,0.36320017261);
\draw(0.940424300357,0.493589088465)--(0.959054704029,0.41906747378);
\draw(0.940424300357,0.493589088465)--(0.940424300357,0.547123155666);
\draw(0.906338263824,0.595847198066)--(0.940424300357,0.493589088465);
\draw(0.906338263824,0.595847198066)--(0.906338263824,0.677941159984);
\draw(0.835969744015,0.736584237683)--(0.906338263824,0.595847198066);
\draw(0.835969744015,0.736584237683)--(0.835969744015,0.875385672435);
\draw(0.670440949108,0.90211303259)--(0.835969744015,0.736584237683);
\draw(0.670440949108,1.16972939267)--(0.715771340427,1.12439900135);
\draw(0.670440949108,1.16972939267)--(0.633478231341,1.2436548282);
\draw(0.670440949108,0.90211303259)--(0.670440949108,1.16972939267);
\draw(0.28407904384,0.90211303259)--(0.670440949108,0.90211303259);
\draw(0.28407904384,1.48494113807)--(0.3779862054,1.48494113807);
\draw(0.200037034439,1.56898314747)--(0.152208853468,1.59289723796);
\draw(0.28407904384,1.48494113807)--(0.200037034439,1.56898314747);
\draw(0.28407904384,0.90211303259)--(0.28407904384,1.48494113807);
\draw(0.,0.61803398875)--(0.28407904384,0.90211303259);
\draw(-0.200037034439,1.56898314747)--(-0.152208853468,1.59289723796);
\draw(-0.28407904384,1.48494113807)--(-0.200037034439,1.56898314747);
\draw(-0.28407904384,1.48494113807)--(-0.3779862054,1.48494113807);
\draw(-0.28407904384,0.90211303259)--(-0.28407904384,1.48494113807);
\draw(-0.670440949108,1.16972939267)--(-0.633478231341,1.2436548282);
\draw(-0.670440949108,1.16972939267)--(-0.715771340427,1.12439900135);
\draw(-0.670440949108,0.90211303259)--(-0.670440949108,1.16972939267);
\draw(-0.835969744015,0.736584237683)--(-0.835969744015,0.875385672435);
\draw(-0.906338263824,0.595847198066)--(-0.906338263824,0.677941159984);
\draw(-0.940424300357,0.493589088465)--(-0.940424300357,0.547123155666);
\draw(-0.959054704029,0.41906747378)--(-0.970228164263,0.36320017261);
\draw(-0.940424300357,0.493589088465)--(-0.959054704029,0.41906747378);
\draw(-0.906338263824,0.595847198066)--(-0.940424300357,0.493589088465);
\draw(-0.835969744015,0.736584237683)--(-0.906338263824,0.595847198066);
\draw(-0.670440949108,0.90211303259)--(-0.835969744015,0.736584237683);
\draw(-0.28407904384,0.90211303259)--(-0.670440949108,0.90211303259);
\draw(0.,0.61803398875)--(-0.28407904384,0.90211303259);
\end{scope}
\end{tikzpicture}}
\caption{Tropical caustic of an ellipse. This particular ellipse is symmetric under sign changes of coordinates, with the ratio between width and height being golden. Removing the central edge, the caustic splits into four congruent caustic branches. The top right branch has the slopes described in \ref{ssec_dulabranch}.}
\label{fig-ellipse}
\end{figure}

\vspace{-5pt}

\subsection{Tropical caustics}
The term ``tropical caustic'', inspired by analogies with Vladimir Arnold's work on geometric optics and symplectic geometry, was proposed by Grigory Mikhalkin and appears in the literature for the first time in \cite{MSthese} as a byproduct of \cite{KS26}. Tropical caustics of planar convex domains are (possibly infinite) planar trees with naturally weighted straight edges of rational slopes -- see Figure \ref{fig-ellipse}, where all edges except for the middle one have weight one, and the middle edge has weight two. In particular, this object may be called a rational tropical (analytic) curve as the balancing condition \cite{BRIEF} is satisfied at every vertex. The adjective ``analytic'' is particularly justified if one uses the most compact definition of a tropical caustic, i.e. as the locus of non-linearity of tropical distance series (in contrast with the algebraic case of a tropical polynomial) constructed from a given convex domain -- this series is given by the infimum (i.e. tropical sum) of all affine support functions (i.e. tropical monomials) of the domain having integer gradients.

\subsection{Branching from the dual subdivision's perspective}
\label{ssec_dulabranch}
Start with the seed triangle having vertices
\[
(0,0),\qquad (0,1),\qquad (1,0).
\]
Declare the zeroth-step edge to be the edge connecting $(1,0)$ and $(0,1)$.

At step $n+1$, for every step-$n$ edge connecting $\lambda_1$ and $\lambda_2$, introduce the new point
\[
\lambda_1+\lambda_2,
\]
and call the two edges
\[
[\lambda_1,\lambda_1+\lambda_2]
\qquad\text{and}\qquad
[\lambda_1+\lambda_2,\lambda_2]
\]
the step-$(n+1)$ edges. There are $2^{n+1}$ such edges.

Every triangle with vertices
\[
\lambda_1,\qquad \lambda_2,\qquad \lambda_1+\lambda_2
\]
has area $\tfrac12$. Call the union of all these triangles the \emph{standard Farey triangulation}. Its nonzero vertices are the primitive vectors in $\Z^2_{\ge 0},$ i.e. such whose coordinates are coprime.% This triangulation is a geometric enhancement of Farey sequence $F_n$, where the members of the sequence are slopes of vertices at step $n$.
This triangulation is related to Farey sequences through the mediant construction. In principle, we could have initiated such a process starting with any seed triangle with a distinguished side.

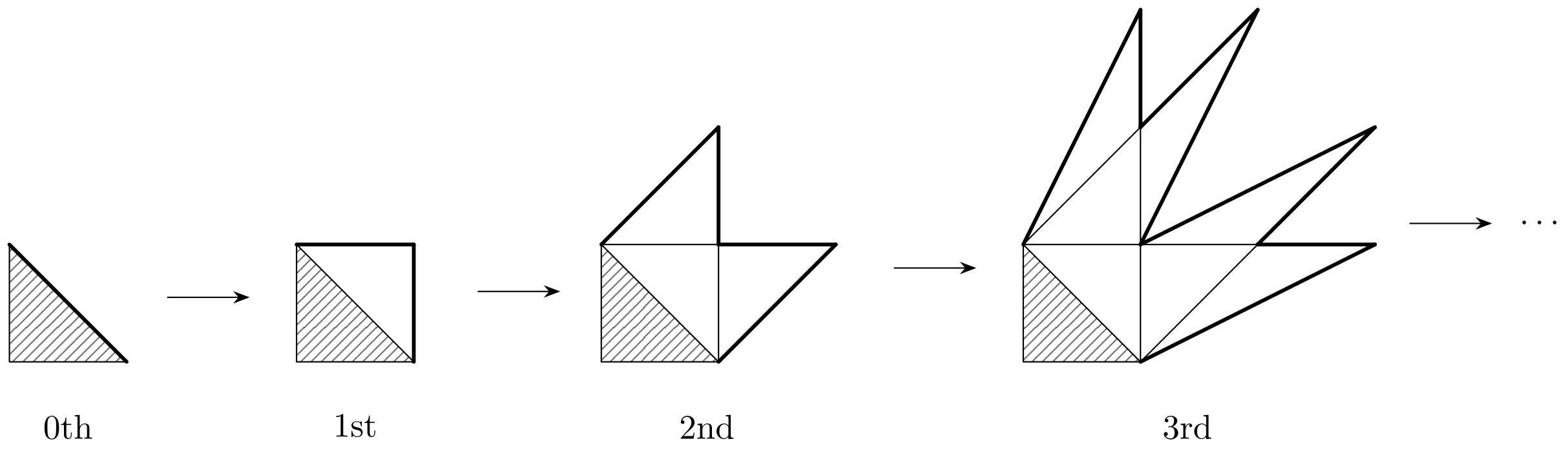
\begin{figure}[H]
\centering
\begin{tikzpicture}[scale=1.2,line cap=round,line join=round,>=Stealth]
  \tikzset{seed/.style={pattern=north east lines,pattern color=black!55}}

  % Stage 0
  \begin{scope}[xshift=0cm]
    \coordinate (O) at (0,0); \coordinate (A) at (1,0); \coordinate (B) at (0,1);
    \fill[seed] (O)--(A)--(B)--cycle;
    \draw (O)--(A)--(B)--cycle;
    \draw[line width=1.1pt] (A)--(B);
    \node[below] at (0.5,-0.35) {$0$th};
  \end{scope}

  \draw[->] (1.35,0.55)--(2.05,0.55);

  % Stage 1
  \begin{scope}[xshift=2.45cm]
    \coordinate (O) at (0,0); \coordinate (A) at (1,0); \coordinate (B) at (0,1); \coordinate (C) at (1,1);
    \fill[seed] (O)--(A)--(B)--cycle;
    \draw (O)--(A)--(B)--cycle;
    \draw (A)--(C)--(B);
    \draw[line width=1.1pt] (A)--(C) (C)--(B);
    \node[below] at (0.5,-0.35) {$1$st};
  \end{scope}

  \draw[->] (4.0,0.6)--(4.7,0.6);

  % Stage 2
  \begin{scope}[xshift=5.05cm]
    \coordinate (O) at (0,0); \coordinate (A) at (1,0); \coordinate (B) at (0,1); \coordinate (C) at (1,1);
    \coordinate (D) at (2,1); \coordinate (E) at (1,2);
    \fill[seed] (O)--(A)--(B)--cycle;
    \draw (O)--(A)--(B)--cycle;
    \draw (A)--(C)--(B);
    \draw (A)--(D)--(C) (C)--(E)--(B);
    \draw[line width=1.1pt] (A)--(D) (D)--(C) (C)--(E) (E)--(B);
    \node[below] at (0.9,-0.35) {$2$nd};
  \end{scope}

  \draw[->] (7.55,0.8)--(8.25,0.8);

  % Stage 3
  \begin{scope}[xshift=8.65cm]
    \coordinate (O) at (0,0); \coordinate (A) at (1,0); \coordinate (B) at (0,1); \coordinate (C) at (1,1);
    \coordinate (D) at (2,1); \coordinate (E) at (1,2);
    \coordinate (F) at (3,1); \coordinate (G) at (3,2); \coordinate (H) at (2,3); \coordinate (I) at (1,3);
    \fill[seed] (O)--(A)--(B)--cycle;
    \draw (O)--(A)--(B)--cycle;
    \draw (A)--(C)--(B);
    \draw (A)--(D)--(C) (C)--(E)--(B);
    \draw (A)--(F)--(D) (D)--(G)--(C) (C)--(H)--(E) (E)--(I)--(B);
    \draw[line width=1.1pt] (A)--(F) (F)--(D) (D)--(G) (G)--(C)
      (C)--(H) (H)--(E) (E)--(I) (I)--(B);
    \node[below] at (1.4,-0.35) {$3$rd};
  \end{scope}

  \draw[->] (11.95,1.18)--(12.65,1.18);
  \node[right] at (12.78,1.18) {$\cdots$};
\end{tikzpicture}
\caption{The first stages of the standard Farey triangulation, its triangles contain all points of the positive quadrant except those which are $r$-multiples of primitive vectors with $r>1$. The seed triangle is drawn shaded.}
\label{fig-FareyTriang}
\end{figure}

\subsection{The dual caustic branch}

The caustic branch is dual to this triangulation, which should be thought of as an analytic analogue of the standard (algebraic) dual subdivision of Newton polygon, see \cite{BRIEF}.

\begin{itemize}
  \item The stem edge of the caustic branch starts at $(0,0)$ and is orthogonal to the zeroth step edge of the Farey triangulation, i.e. it goes in the direction $(1,1).$
  \item Pass to the first triangle and branch toward the other two sides of that triangle.
  \item At step $2$, the two new caustic edges branch orthogonally to the corresponding triangulation edges; continue in the same manner at every subsequent step.
\end{itemize}

Thus, the $1$-skeleton of the triangulation is dual to the caustic branching, as indicated schematically below. Tropical caustic curves are made of such branches, see \cite{MIKH}.

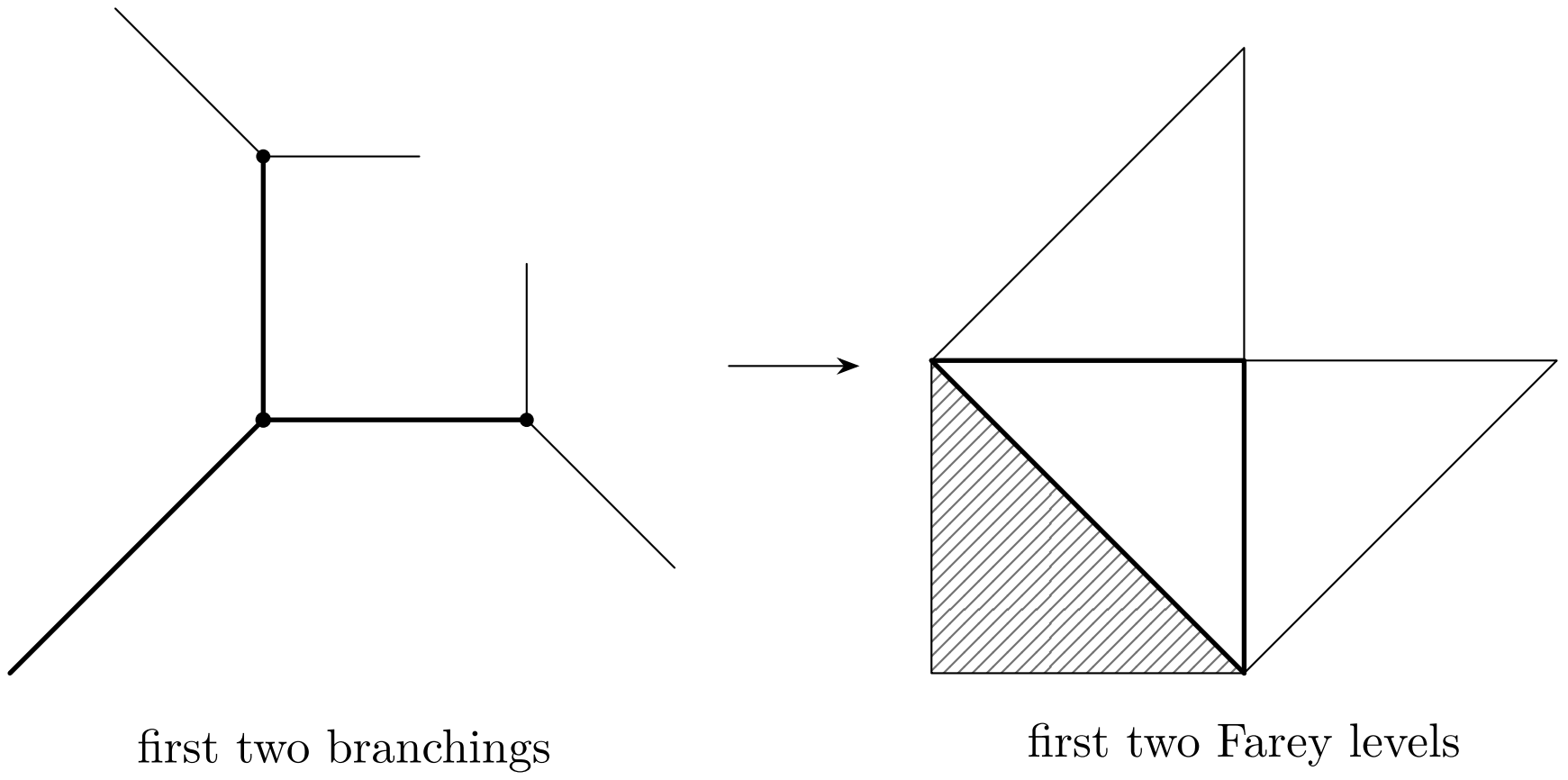
\begin{figure}[H]
\centering
\begin{tikzpicture}[scale=0.96,>=Stealth,line cap=round,line join=round]
  \tikzset{seed/.style={pattern=north east lines,pattern color=black!55}}
  % larger caustic branch, no labels on edges
  \begin{scope}[xshift=-0.20cm,yshift=0.00cm,scale=1.55]
    \coordinate (O) at (0,0);
    \coordinate (V) at (1.25,1.25);
    \coordinate (U) at (1.25,2.55);
    \coordinate (R) at (2.55,1.25);
    \coordinate (UU) at (0.52,3.28);
    \coordinate (UR) at (2.02,2.55);
    \coordinate (RU) at (2.55,2.02);
    \coordinate (RR) at (3.28,0.52);

    \draw[line width=1pt] (O)--(V);
    \draw[line width=1pt] (V)--(U) (V)--(R);
    \draw[thin] (U)--(UU) (U)--(UR);
    \draw[thin] (R)--(RU) (R)--(RR);
    \fill (V) circle (1.1pt);
    \fill (U) circle (1.0pt);
    \fill (R) circle (1.0pt);
    \node[below=4pt] at (1.65,-0.10) {first two branchings};
  \end{scope}

  \draw[->] (5.30,2.35)--(6.30,2.35);

  % triangulation with four visible triangles and hatched seed triangle
  \begin{scope}[xshift=6.85cm,yshift=0cm,scale=0.92]
    \coordinate (A) at (0,0);
    \coordinate (B) at (0,2.6);
    \coordinate (C) at (2.6,0);
    \coordinate (D) at (2.6,2.6);
    \coordinate (E) at (5.2,2.6);
    \coordinate (F) at (2.6,5.2);

    \fill[seed] (A)--(B)--(C)--cycle;
    \draw (A)--(B)--(C)--cycle;
    \draw[line width=1pt] (B)--(C);
    \draw[line width=1pt] (B)--(D)--(C);
    \draw (B)--(F)--(D);
    \draw (D)--(E)--(C);
    \node[below=4pt] at (2.6,-0.12) {first two Farey levels};
  \end{scope}
\end{tikzpicture}
\caption{The caustic branch determines the first two Farey levels.}
\end{figure}

\subsection{The length rule}

In the construction above, the lengths of the caustic edges were implicit. The rule is as follows.

\begin{itemize}
  \item  Define tropical length of primitive vectors to be $1$.
  \item For every caustic vertex $v$ dual to the triangle with vertices
  $\lambda_1$, $\lambda_2$, and $\lambda_1+\lambda_2$, assign a positive number
  \[
  a_v=\potential{\lambda_1}{\lambda_2}>0,
  \]
  called its \emph{potential}, such that $a_{v_n}\longrightarrow 0$ over any sequence $\{v_n\}_{n=1}^{\infty},$ where $v_n$ is a vertex of generation $n.$   
  \item For a caustic edge whose endpoints $v_{n}$ and $v_{n+1}$ arise at steps $n$ and $n+1$, respectively, take its tropical length to be $a_{v_{n}}-a_{v_{n+1}}.$ Thus the total length of any complete geodesic on the branch is $a_{(0,0)}.$ 
\end{itemize}
\begin{definition}
    For compact convex domain $\Omega\subset\mathbb{R}^2,$ its tropical distance series is defined as $$\mathcal{F}_\Omega(z)=\inf_{\lambda\in\mathbb{Z}^2\setminus\{(0,0)\}}(c_\lambda+\lambda\cdot z),$$ where $c_\lambda=-\inf_{z\in\Omega}(\lambda\cdot z)$ and $z \in \mathbb{R}^2$.  Note that $\Omega$ is the domain of convergence of $\mathcal{F}_\Omega,$ i.e. the value $\mathcal{F}_\Omega(z)$ is finite if and only if $z\in\Omega.$ Tropical caustic of $\Omega$ is the set of points $\mathcal{K}_\Omega$ in $\Omega$ where $\mathcal{F}_\Omega$ is not locally affine. Tropical caustic of a compact convex domain is always a tree, and if $\Omega$ has smooth boundary, then removing from $\mathcal{K}_\Omega$ the locus where $\mathcal{F}_\Omega$ attains its maximum one gets a union of branches all of which are $\operatorname{SL}_2(\mathbb{Z})$ and translation congruent to the one described above. The potential $a_v$ at vertex of the caustic $v$ is precisely $\mathcal{F}_\Omega(v).$
\end{definition}

\subsection{Examples}

\subsubsection*{Example 1: Caustic branch of unit disc}

Take
\begin{equation}
\potential{\lambda_1}{\lambda_2}
=
\norm{\lambda_1}_2+
\norm{\lambda_2}_2-
\norm{\lambda_1+\lambda_2}_2,
\label{eq:euclidean-potential}
\end{equation}
where $\norm{\cdot}_2$ is the Euclidean norm. It appears, see \cite{KSpi}, that the limit set of the corresponding caustic branch is
\[
\left\{\left(x,\sqrt{1-x^2}\right):x\in[0,1]\right\},
\]
the quarter-circle arc, see the left part of Figure \ref{fig7}.

\begin{remark}
    More generally, replacing the Euclidean norm with the square root of a positive definite quadratic form $Q(x,y)$, one gets an arc of $Q^*=1$ as the limit set, where $Q^*$ is the dual quadratic form to $Q$ --- combining with determining the minimal model from symmetries, this is the way how Figure \ref{fig-ellipse} was calculated.
\end{remark}%one gets не ясно в рамках какого действия, до этого ты вроде не говорил что нас интересуют предельные точки

\begin{remark}
    Similarly, for the replacement of the Euclidean $2$-norm by a $p$-norm with $p\geq 1$ we get the positive arc of $\norm{(x,y)}_q=1$ with the usual relation $p^{-1}+q^{-1}=1$.
\end{remark}

\subsubsection*{Example 2: Caustic branch of domain L}

Take
\begin{equation}
\potential{\lambda_1}{\lambda_2}
=
\frac{1}{
\norm{\lambda_1}_1\,
\cdot \norm{\lambda_2}_1\,
\cdot \norm{\lambda_1+\lambda_2}_1},
\label{eq:l1-potential}
\end{equation}
where
\[
\norm{(p,q)}_1=|p|+|q|.
\]
Then the set of the limiting points is (see \cite{KSsummation} for the proof) as shown in the right part of Figure \ref{fig7}:
\[
\left\{(x,y):\sqrt{1-x}+\sqrt{1-y}=1,\quad x,y\in[0,1]\right\}.
\]
\begin{remark}
    In particular, this is the limit shape of random lattice polygons in a square when the mesh of the lattice goes to zero: Vershik, B\'ar\'any and Sinai wrote on this and related properties in 1994, but the rescaled curve itself belongs to the class of curves maximizing the number of lattice points on it, this problem was studied by Jarn\'ik already in 1926. So far, this is the only smooth curve for which we know the pole to the left from $\frac{2}{3}$ of its tropical zeta function, see section 5 of \cite{Zeta}, where a few more geometrical and physical references are mentioned.
\end{remark}

\begin{figure}[H]
\centering
\begin{tikzpicture}[>=Stealth,line cap=round,line join=round]
  \begin{scope}[xshift=-4.35cm,yshift=0cm,scale=2.72]
    \draw[->] (0,0)--(1.12,0) node[right] {$x$};
    \draw[->] (0,0)--(0,1.12) node[above] {$y$};
    \draw (1,-0.02)--(1,0.02) node[below=5pt] {$1$};
    \draw (-0.02,1)--(0.02,1) node[left=5pt] {$1$};
    \draw[line width=1pt,samples=140,domain=0:90,smooth]
      plot ({cos(\x)},{sin(\x)});
    \draw[line width=0.95pt] (0.00000000,0.00000000)--(0.41421356,0.41421356);
    \draw[line width=0.95pt] (0.41421356,0.41421356)--(0.82185442,0.41421356);
    \draw[line width=0.62pt] (0.82185442,0.41421356)--(0.92620968,0.30985829);
    \draw[line width=0.62pt] (0.92620968,0.30985829)--(0.96082797,0.24062173);
    \draw[line width=0.43pt] (0.96082797,0.24062173)--(0.97591389,0.19536396);
    \draw[line width=0.43pt] (0.97591389,0.19536396)--(0.98374302,0.16404745);
    \draw[line width=0.43pt] (0.98374302,0.16404745)--(0.98830528,0.14123612);
    \draw[line width=0.43pt] (0.98374302,0.16404745)--(0.98374302,0.17886227);
    \draw[line width=0.43pt] (0.97591389,0.19536396)--(0.97591389,0.21686939);
    \draw[line width=0.43pt] (0.97591389,0.21686939)--(0.97775166,0.20951829);
    \draw[line width=0.43pt] (0.97591389,0.21686939)--(0.97431922,0.22484271);
    \draw[line width=0.43pt] (0.96082797,0.24062173)--(0.96082797,0.27452037);
    \draw[line width=0.43pt] (0.96082797,0.27452037)--(0.96464010,0.26308397);
    \draw[line width=0.43pt] (0.96464010,0.26308397)--(0.96619227,0.25765138);
    \draw[line width=0.43pt] (0.96464010,0.26308397)--(0.96350128,0.26763923);
    \draw[line width=0.43pt] (0.96082797,0.27452037)--(0.95763561,0.28728979);
    \draw[line width=0.43pt] (0.95763561,0.28728979)--(0.95934540,0.28216042);
    \draw[line width=0.43pt] (0.95763561,0.28728979)--(0.95570085,0.29406145);
    \draw[line width=0.62pt] (0.92620968,0.30985829)--(0.92620968,0.37046778);
    \draw[line width=0.43pt] (0.92620968,0.37046778)--(0.93595179,0.35098356);
    \draw[line width=0.43pt] (0.93595179,0.35098356)--(0.93966625,0.34169743);
    \draw[line width=0.43pt] (0.93966625,0.34169743)--(0.94168426,0.33631606);
    \draw[line width=0.43pt] (0.93966625,0.34169743)--(0.93833168,0.34570114);
    \draw[line width=0.43pt] (0.93595179,0.35098356)--(0.93329334,0.35895890);
    \draw[line width=0.43pt] (0.93329334,0.35895890)--(0.93447996,0.35599237);
    \draw[line width=0.43pt] (0.93329334,0.35895890)--(0.93198802,0.36243976);
    \draw[line width=0.43pt] (0.92620968,0.37046778)--(0.91848853,0.39363124);
    \draw[line width=0.43pt] (0.91848853,0.39363124)--(0.92301030,0.38458771);
    \draw[line width=0.43pt] (0.92301030,0.38458771)--(0.92465855,0.38074179);
    \draw[line width=0.43pt] (0.92301030,0.38458771)--(0.92162443,0.38805238);
    \draw[line width=0.43pt] (0.91848853,0.39363124)--(0.91353574,0.40601322);
    \draw[line width=0.43pt] (0.91353574,0.40601322)--(0.91613731,0.40081008);
    \draw[line width=0.43pt] (0.91353574,0.40601322)--(0.91022531,0.41373755);
    \draw[line width=0.62pt] (0.82185442,0.41421356)--(0.82185442,0.54762888);
    \draw[line width=0.62pt] (0.82185442,0.54762888)--(0.85591732,0.51356598);
    \draw[line width=0.43pt] (0.85591732,0.51356598)--(0.86772779,0.49585027);
    \draw[line width=0.43pt] (0.86772779,0.49585027)--(0.87392741,0.48551758);
    \draw[line width=0.43pt] (0.87392741,0.48551758)--(0.87777362,0.47878671);
    \draw[line width=0.43pt] (0.87392741,0.48551758)--(0.87155996,0.49025247);
    \draw[line width=0.43pt] (0.86772779,0.49585027)--(0.86373132,0.50384322);
    \draw[line width=0.43pt] (0.86373132,0.50384322)--(0.86541827,0.50103163);
    \draw[line width=0.43pt] (0.86373132,0.50384322)--(0.86192055,0.50701206);
    \draw[line width=0.43pt] (0.85591732,0.51356598)--(0.84777200,0.52985662);
    \draw[line width=0.43pt] (0.84777200,0.52985662)--(0.85162507,0.52407701);
    \draw[line width=0.43pt] (0.85162507,0.52407701)--(0.85327105,0.52144344);
    \draw[line width=0.43pt] (0.85162507,0.52407701)--(0.85026031,0.52635162);
    \draw[line width=0.43pt] (0.84777200,0.52985662)--(0.84358868,0.53682882);
    \draw[line width=0.43pt] (0.84358868,0.53682882)--(0.84548210,0.53398870);
    \draw[line width=0.43pt] (0.84358868,0.53682882)--(0.84114624,0.54073673);
    \draw[line width=0.62pt] (0.82185442,0.54762888)--(0.79688899,0.59755974);
    \draw[line width=0.43pt] (0.79688899,0.59755974)--(0.81342782,0.58102091);
    \draw[line width=0.43pt] (0.81342782,0.58102091)--(0.81914306,0.57340058);
    \draw[line width=0.43pt] (0.81914306,0.57340058)--(0.82215475,0.56918423);
    \draw[line width=0.43pt] (0.81914306,0.57340058)--(0.81695888,0.57667685);
    \draw[line width=0.43pt] (0.81342782,0.58102091)--(0.80868532,0.58813466);
    \draw[line width=0.43pt] (0.80868532,0.58813466)--(0.81067220,0.58548548);
    \draw[line width=0.43pt] (0.80868532,0.58813466)--(0.80638788,0.59135107);
    \draw[line width=0.43pt] (0.79688899,0.59755974)--(0.77953796,0.62358627);
    \draw[line width=0.43pt] (0.77953796,0.62358627)--(0.78925730,0.61386694);
    \draw[line width=0.43pt] (0.78925730,0.61386694)--(0.79259811,0.60969092);
    \draw[line width=0.43pt] (0.78925730,0.61386694)--(0.78630136,0.61780819);
    \draw[line width=0.43pt] (0.77953796,0.62358627)--(0.76753436,0.63959108);
    \draw[line width=0.43pt] (0.76753436,0.63959108)--(0.77391897,0.63320646);
    \draw[line width=0.43pt] (0.76753436,0.63959108)--(0.75885699,0.65043779);
    \draw[line width=0.95pt] (0.41421356,0.41421356)--(0.41421356,0.82185442);
    \draw[line width=0.62pt] (0.41421356,0.82185442)--(0.54762888,0.82185442);
    \draw[line width=0.62pt] (0.54762888,0.82185442)--(0.59755974,0.79688899);
    \draw[line width=0.43pt] (0.59755974,0.79688899)--(0.62358627,0.77953796);
    \draw[line width=0.43pt] (0.62358627,0.77953796)--(0.63959108,0.76753436);
    \draw[line width=0.43pt] (0.63959108,0.76753436)--(0.65043779,0.75885699);
    \draw[line width=0.43pt] (0.63959108,0.76753436)--(0.63320646,0.77391897);
    \draw[line width=0.43pt] (0.62358627,0.77953796)--(0.61386694,0.78925730);
    \draw[line width=0.43pt] (0.61386694,0.78925730)--(0.61780819,0.78630136);
    \draw[line width=0.43pt] (0.61386694,0.78925730)--(0.60969092,0.79259811);
    \draw[line width=0.43pt] (0.59755974,0.79688899)--(0.58102091,0.81342782);
    \draw[line width=0.43pt] (0.58102091,0.81342782)--(0.58813466,0.80868532);
    \draw[line width=0.43pt] (0.58813466,0.80868532)--(0.59135107,0.80638788);
    \draw[line width=0.43pt] (0.58813466,0.80868532)--(0.58548548,0.81067220);
    \draw[line width=0.43pt] (0.58102091,0.81342782)--(0.57340058,0.81914306);
    \draw[line width=0.43pt] (0.57340058,0.81914306)--(0.57667685,0.81695888);
    \draw[line width=0.43pt] (0.57340058,0.81914306)--(0.56918423,0.82215475);
    \draw[line width=0.62pt] (0.54762888,0.82185442)--(0.51356598,0.85591732);
    \draw[line width=0.43pt] (0.51356598,0.85591732)--(0.52985662,0.84777200);
    \draw[line width=0.43pt] (0.52985662,0.84777200)--(0.53682882,0.84358868);
    \draw[line width=0.43pt] (0.53682882,0.84358868)--(0.54073673,0.84114624);
    \draw[line width=0.43pt] (0.53682882,0.84358868)--(0.53398870,0.84548210);
    \draw[line width=0.43pt] (0.52985662,0.84777200)--(0.52407701,0.85162507);
    \draw[line width=0.43pt] (0.52407701,0.85162507)--(0.52635162,0.85026031);
    \draw[line width=0.43pt] (0.52407701,0.85162507)--(0.52144344,0.85327105);
    \draw[line width=0.43pt] (0.51356598,0.85591732)--(0.49585027,0.86772779);
    \draw[line width=0.43pt] (0.49585027,0.86772779)--(0.50384322,0.86373132);
    \draw[line width=0.43pt] (0.50384322,0.86373132)--(0.50701206,0.86192055);
    \draw[line width=0.43pt] (0.50384322,0.86373132)--(0.50103163,0.86541827);
    \draw[line width=0.43pt] (0.49585027,0.86772779)--(0.48551758,0.87392741);
    \draw[line width=0.43pt] (0.48551758,0.87392741)--(0.49025247,0.87155996);
    \draw[line width=0.43pt] (0.48551758,0.87392741)--(0.47878671,0.87777362);
    \draw[line width=0.62pt] (0.41421356,0.82185442)--(0.30985829,0.92620968);
    \draw[line width=0.62pt] (0.30985829,0.92620968)--(0.37046778,0.92620968);
    \draw[line width=0.43pt] (0.37046778,0.92620968)--(0.39363124,0.91848853);
    \draw[line width=0.43pt] (0.39363124,0.91848853)--(0.40601322,0.91353574);
    \draw[line width=0.43pt] (0.40601322,0.91353574)--(0.41373755,0.91022531);
    \draw[line width=0.43pt] (0.40601322,0.91353574)--(0.40081008,0.91613731);
    \draw[line width=0.43pt] (0.39363124,0.91848853)--(0.38458771,0.92301030);
    \draw[line width=0.43pt] (0.38458771,0.92301030)--(0.38805238,0.92162443);
    \draw[line width=0.43pt] (0.38458771,0.92301030)--(0.38074179,0.92465855);
    \draw[line width=0.43pt] (0.37046778,0.92620968)--(0.35098356,0.93595179);
    \draw[line width=0.43pt] (0.35098356,0.93595179)--(0.35895890,0.93329334);
    \draw[line width=0.43pt] (0.35895890,0.93329334)--(0.36243976,0.93198802);
    \draw[line width=0.43pt] (0.35895890,0.93329334)--(0.35599237,0.93447996);
    \draw[line width=0.43pt] (0.35098356,0.93595179)--(0.34169743,0.93966625);
    \draw[line width=0.43pt] (0.34169743,0.93966625)--(0.34570114,0.93833168);
    \draw[line width=0.43pt] (0.34169743,0.93966625)--(0.33631606,0.94168426);
    \draw[line width=0.62pt] (0.30985829,0.92620968)--(0.24062173,0.96082797);
    \draw[line width=0.43pt] (0.24062173,0.96082797)--(0.27452037,0.96082797);
    \draw[line width=0.43pt] (0.27452037,0.96082797)--(0.28728979,0.95763561);
    \draw[line width=0.43pt] (0.28728979,0.95763561)--(0.29406145,0.95570085);
    \draw[line width=0.43pt] (0.28728979,0.95763561)--(0.28216042,0.95934540);
    \draw[line width=0.43pt] (0.27452037,0.96082797)--(0.26308397,0.96464010);
    \draw[line width=0.43pt] (0.26308397,0.96464010)--(0.26763923,0.96350128);
    \draw[line width=0.43pt] (0.26308397,0.96464010)--(0.25765138,0.96619227);
    \draw[line width=0.43pt] (0.24062173,0.96082797)--(0.19536396,0.97591389);
    \draw[line width=0.43pt] (0.19536396,0.97591389)--(0.21686939,0.97591389);
    \draw[line width=0.43pt] (0.21686939,0.97591389)--(0.22484271,0.97431922);
    \draw[line width=0.43pt] (0.21686939,0.97591389)--(0.20951829,0.97775166);
    \draw[line width=0.43pt] (0.19536396,0.97591389)--(0.16404745,0.98374302);
    \draw[line width=0.43pt] (0.16404745,0.98374302)--(0.17886227,0.98374302);
    \draw[line width=0.43pt] (0.16404745,0.98374302)--(0.14123612,0.98830528);
    \node[below=11pt] at (0.50,0) {$y=\sqrt{1-x^2},\ x\geq 0$};
  \end{scope}

  \begin{scope}[xshift=4.35cm,yshift=0cm,scale=2.72]
    \draw[->] (0,0)--(1.12,0) node[right] {$x$};
    \draw[->] (0,0)--(0,1.12) node[above] {$y$};
    \draw (1,-0.02)--(1,0.02) node[below=5pt] {$1$};
    \draw (-0.02,1)--(0.02,1) node[left=5pt] {$1$};
    \draw[line width=1pt,samples=180,domain=0:1,smooth]
      plot ({\x},{\x-1+2*sqrt(1-\x)});
    \draw[line width=0.95pt] (0.00000000,0.00000000)--(0.50000000,0.50000000);
    \draw[line width=0.95pt] (0.50000000,0.50000000)--(0.83333333,0.50000000);
    \draw[line width=0.62pt] (0.83333333,0.50000000)--(0.91666667,0.41666667);
    \draw[line width=0.62pt] (0.91666667,0.41666667)--(0.95000000,0.35000000);
    \draw[line width=0.43pt] (0.95000000,0.35000000)--(0.96666667,0.30000000);
    \draw[line width=0.43pt] (0.96666667,0.30000000)--(0.97619048,0.26190476);
    \draw[line width=0.43pt] (0.97619048,0.26190476)--(0.98214286,0.23214286);
    \draw[line width=0.43pt] (0.97619048,0.26190476)--(0.97619048,0.28388278);
    \draw[line width=0.43pt] (0.96666667,0.30000000)--(0.96666667,0.33030303);
    \draw[line width=0.43pt] (0.96666667,0.33030303)--(0.96880570,0.32174688);
    \draw[line width=0.43pt] (0.96666667,0.33030303)--(0.96477273,0.33977273);
    \draw[line width=0.43pt] (0.95000000,0.35000000)--(0.95000000,0.39444444);
    \draw[line width=0.43pt] (0.95000000,0.39444444)--(0.95396825,0.38253968);
    \draw[line width=0.43pt] (0.95396825,0.38253968)--(0.95563910,0.37669173);
    \draw[line width=0.43pt] (0.95396825,0.38253968)--(0.95272602,0.38750863);
    \draw[line width=0.43pt] (0.95000000,0.39444444)--(0.94658120,0.40811966);
    \draw[line width=0.43pt] (0.94658120,0.40811966)--(0.94832945,0.40287490);
    \draw[line width=0.43pt] (0.94658120,0.40811966)--(0.94457014,0.41515837);
    \draw[line width=0.62pt] (0.91666667,0.41666667)--(0.91666667,0.48809524);
    \draw[line width=0.43pt] (0.91666667,0.48809524)--(0.92532468,0.47077922);
    \draw[line width=0.43pt] (0.92532468,0.47077922)--(0.92878788,0.46212121);
    \draw[line width=0.43pt] (0.92878788,0.46212121)--(0.93070175,0.45701754);
    \draw[line width=0.43pt] (0.92878788,0.46212121)--(0.92750583,0.46596737);
    \draw[line width=0.43pt] (0.92532468,0.47077922)--(0.92279942,0.47835498);
    \draw[line width=0.43pt] (0.92279942,0.47835498)--(0.92389411,0.47561825);
    \draw[line width=0.43pt] (0.92279942,0.47835498)--(0.92158730,0.48158730);
    \draw[line width=0.43pt] (0.91666667,0.48809524)--(0.90952381,0.50952381);
    \draw[line width=0.43pt] (0.90952381,0.50952381)--(0.91344538,0.50168067);
    \draw[line width=0.43pt] (0.91344538,0.50168067)--(0.91491597,0.49824930);
    \draw[line width=0.43pt] (0.91344538,0.50168067)--(0.91220044,0.50479303);
    \draw[line width=0.43pt] (0.90952381,0.50952381)--(0.90512821,0.52051282);
    \draw[line width=0.43pt] (0.90512821,0.52051282)--(0.90735786,0.51605351);
    \draw[line width=0.43pt] (0.90512821,0.52051282)--(0.90224359,0.52724359);
    \draw[line width=0.62pt] (0.83333333,0.50000000)--(0.83333333,0.63333333);
    \draw[line width=0.62pt] (0.83333333,0.63333333)--(0.85833333,0.60833333);
    \draw[line width=0.43pt] (0.85833333,0.60833333)--(0.86742424,0.59469697);
    \draw[line width=0.43pt] (0.86742424,0.59469697)--(0.87229437,0.58658009);
    \draw[line width=0.43pt] (0.87229437,0.58658009)--(0.87535014,0.58123249);
    \draw[line width=0.43pt] (0.87229437,0.58658009)--(0.87038961,0.59038961);
    \draw[line width=0.43pt] (0.86742424,0.59469697)--(0.86423445,0.60107656);
    \draw[line width=0.43pt] (0.86423445,0.60107656)--(0.86555024,0.59888357);
    \draw[line width=0.43pt] (0.86423445,0.60107656)--(0.86281676,0.60355750);
    \draw[line width=0.43pt] (0.85833333,0.60833333)--(0.85192308,0.62115385);
    \draw[line width=0.43pt] (0.85192308,0.62115385)--(0.85485348,0.61675824);
    \draw[line width=0.43pt] (0.85485348,0.61675824)--(0.85611658,0.61473727);
    \draw[line width=0.43pt] (0.85485348,0.61675824)--(0.85380306,0.61850894);
    \draw[line width=0.43pt] (0.85192308,0.62115385)--(0.84871795,0.62649573);
    \draw[line width=0.43pt] (0.84871795,0.62649573)--(0.85015164,0.62434519);
    \draw[line width=0.43pt] (0.84871795,0.62649573)--(0.84685990,0.62946860);
    \draw[line width=0.62pt] (0.83333333,0.63333333)--(0.81428571,0.67142857);
    \draw[line width=0.43pt] (0.81428571,0.67142857)--(0.82619048,0.65952381);
    \draw[line width=0.43pt] (0.82619048,0.65952381)--(0.83039216,0.65392157);
    \draw[line width=0.43pt] (0.83039216,0.65392157)--(0.83262032,0.65080214);
    \draw[line width=0.43pt] (0.83039216,0.65392157)--(0.82876944,0.65635565);
    \draw[line width=0.43pt] (0.82619048,0.65952381)--(0.82268170,0.66478697);
    \draw[line width=0.43pt] (0.82268170,0.66478697)--(0.82413696,0.66284663);
    \draw[line width=0.43pt] (0.82268170,0.66478697)--(0.82099479,0.66714864);
    \draw[line width=0.43pt] (0.81428571,0.67142857)--(0.80158730,0.69047619);
    \draw[line width=0.43pt] (0.80158730,0.69047619)--(0.80853175,0.68353175);
    \draw[line width=0.43pt] (0.80853175,0.68353175)--(0.81094720,0.68051242);
    \draw[line width=0.43pt] (0.80853175,0.68353175)--(0.80638889,0.68638889);
    \draw[line width=0.43pt] (0.80158730,0.69047619)--(0.79292929,0.70202020);
    \draw[line width=0.43pt] (0.79292929,0.70202020)--(0.79747475,0.69747475);
    \draw[line width=0.43pt] (0.79292929,0.70202020)--(0.78671329,0.70979021);
    \draw[line width=0.95pt] (0.50000000,0.50000000)--(0.50000000,0.83333333);
    \draw[line width=0.62pt] (0.50000000,0.83333333)--(0.63333333,0.83333333);
    \draw[line width=0.62pt] (0.63333333,0.83333333)--(0.67142857,0.81428571);
    \draw[line width=0.43pt] (0.67142857,0.81428571)--(0.69047619,0.80158730);
    \draw[line width=0.43pt] (0.69047619,0.80158730)--(0.70202020,0.79292929);
    \draw[line width=0.43pt] (0.70202020,0.79292929)--(0.70979021,0.78671329);
    \draw[line width=0.43pt] (0.70202020,0.79292929)--(0.69747475,0.79747475);
    \draw[line width=0.43pt] (0.69047619,0.80158730)--(0.68353175,0.80853175);
    \draw[line width=0.43pt] (0.68353175,0.80853175)--(0.68638889,0.80638889);
    \draw[line width=0.43pt] (0.68353175,0.80853175)--(0.68051242,0.81094720);
    \draw[line width=0.43pt] (0.67142857,0.81428571)--(0.65952381,0.82619048);
    \draw[line width=0.43pt] (0.65952381,0.82619048)--(0.66478697,0.82268170);
    \draw[line width=0.43pt] (0.66478697,0.82268170)--(0.66714864,0.82099479);
    \draw[line width=0.43pt] (0.66478697,0.82268170)--(0.66284663,0.82413696);
    \draw[line width=0.43pt] (0.65952381,0.82619048)--(0.65392157,0.83039216);
    \draw[line width=0.43pt] (0.65392157,0.83039216)--(0.65635565,0.82876944);
    \draw[line width=0.43pt] (0.65392157,0.83039216)--(0.65080214,0.83262032);
    \draw[line width=0.62pt] (0.63333333,0.83333333)--(0.60833333,0.85833333);
    \draw[line width=0.43pt] (0.60833333,0.85833333)--(0.62115385,0.85192308);
    \draw[line width=0.43pt] (0.62115385,0.85192308)--(0.62649573,0.84871795);
    \draw[line width=0.43pt] (0.62649573,0.84871795)--(0.62946860,0.84685990);
    \draw[line width=0.43pt] (0.62649573,0.84871795)--(0.62434519,0.85015164);
    \draw[line width=0.43pt] (0.62115385,0.85192308)--(0.61675824,0.85485348);
    \draw[line width=0.43pt] (0.61675824,0.85485348)--(0.61850894,0.85380306);
    \draw[line width=0.43pt] (0.61675824,0.85485348)--(0.61473727,0.85611658);
    \draw[line width=0.43pt] (0.60833333,0.85833333)--(0.59469697,0.86742424);
    \draw[line width=0.43pt] (0.59469697,0.86742424)--(0.60107656,0.86423445);
    \draw[line width=0.43pt] (0.60107656,0.86423445)--(0.60355750,0.86281676);
    \draw[line width=0.43pt] (0.60107656,0.86423445)--(0.59888357,0.86555024);
    \draw[line width=0.43pt] (0.59469697,0.86742424)--(0.58658009,0.87229437);
    \draw[line width=0.43pt] (0.58658009,0.87229437)--(0.59038961,0.87038961);
    \draw[line width=0.43pt] (0.58658009,0.87229437)--(0.58123249,0.87535014);
    \draw[line width=0.62pt] (0.50000000,0.83333333)--(0.41666667,0.91666667);
    \draw[line width=0.62pt] (0.41666667,0.91666667)--(0.48809524,0.91666667);
    \draw[line width=0.43pt] (0.48809524,0.91666667)--(0.50952381,0.90952381);
    \draw[line width=0.43pt] (0.50952381,0.90952381)--(0.52051282,0.90512821);
    \draw[line width=0.43pt] (0.52051282,0.90512821)--(0.52724359,0.90224359);
    \draw[line width=0.43pt] (0.52051282,0.90512821)--(0.51605351,0.90735786);
    \draw[line width=0.43pt] (0.50952381,0.90952381)--(0.50168067,0.91344538);
    \draw[line width=0.43pt] (0.50168067,0.91344538)--(0.50479303,0.91220044);
    \draw[line width=0.43pt] (0.50168067,0.91344538)--(0.49824930,0.91491597);
    \draw[line width=0.43pt] (0.48809524,0.91666667)--(0.47077922,0.92532468);
    \draw[line width=0.43pt] (0.47077922,0.92532468)--(0.47835498,0.92279942);
    \draw[line width=0.43pt] (0.47835498,0.92279942)--(0.48158730,0.92158730);
    \draw[line width=0.43pt] (0.47835498,0.92279942)--(0.47561825,0.92389411);
    \draw[line width=0.43pt] (0.47077922,0.92532468)--(0.46212121,0.92878788);
    \draw[line width=0.43pt] (0.46212121,0.92878788)--(0.46596737,0.92750583);
    \draw[line width=0.43pt] (0.46212121,0.92878788)--(0.45701754,0.93070175);
    \draw[line width=0.62pt] (0.41666667,0.91666667)--(0.35000000,0.95000000);
    \draw[line width=0.43pt] (0.35000000,0.95000000)--(0.39444444,0.95000000);
    \draw[line width=0.43pt] (0.39444444,0.95000000)--(0.40811966,0.94658120);
    \draw[line width=0.43pt] (0.40811966,0.94658120)--(0.41515837,0.94457014);
    \draw[line width=0.43pt] (0.40811966,0.94658120)--(0.40287490,0.94832945);
    \draw[line width=0.43pt] (0.39444444,0.95000000)--(0.38253968,0.95396825);
    \draw[line width=0.43pt] (0.38253968,0.95396825)--(0.38750863,0.95272602);
    \draw[line width=0.43pt] (0.38253968,0.95396825)--(0.37669173,0.95563910);
    \draw[line width=0.43pt] (0.35000000,0.95000000)--(0.30000000,0.96666667);
    \draw[line width=0.43pt] (0.30000000,0.96666667)--(0.33030303,0.96666667);
    \draw[line width=0.43pt] (0.33030303,0.96666667)--(0.33977273,0.96477273);
    \draw[line width=0.43pt] (0.33030303,0.96666667)--(0.32174688,0.96880570);
    \draw[line width=0.43pt] (0.30000000,0.96666667)--(0.26190476,0.97619048);
    \draw[line width=0.43pt] (0.26190476,0.97619048)--(0.28388278,0.97619048);
    \draw[line width=0.43pt] (0.26190476,0.97619048)--(0.23214286,0.98214286);
    \node[below=11pt] at (0.50,0) {$\sqrt{1-x}+\sqrt{1-y}=1$};
  \end{scope}

\end{tikzpicture}
\caption{Finite caustic branches computed from the potentials in formulas \eqref{eq:euclidean-potential} and \eqref{eq:l1-potential}. Six Farey levels are shown; the vertices are obtained from the corresponding support lines, so every infinite branch converges to the stated limiting arc.}
  \label{fig7}
\end{figure}

\subsection{The convexity inequalities}

In the current formulation, one assumes that the potentials are monotone along adjacent vertices:
\begin{equation}
   \potential{\lambda_1}{\lambda_2}\geq \potential{\lambda_1+\lambda_2}{\lambda_2}, \qquad\potential{\lambda_1}{\lambda_2}\geq \potential{\lambda_1}{\lambda_1+\lambda_2}.
\end{equation}
  This corresponds to caustic edges going in a ``positive'' direction. In terms of the balancing condition, for monotone potentials all the weights of edges are equal to one, and whenever the monotonicity is violated and the edge goes backwards, we choose the corresponding edge weight to be minus one, resulting in an object which is not a tropical curve in the strict sense, but a local tropical cycle, i.e., since the balancing is still preserved, the local intersection numbers are invariant with respect to deformations.

Crucially, the other inequalities currently hard-wired into tropical optics of convex domains are
\begin{equation}
\potential{\lambda_1}{\lambda_2}
\ge
\sum_{n=0}^{\infty}
\left(
\potential{\lambda_1+\lambda_2}{\lambda_2+n(\lambda_1+\lambda_2)}
+
\potential{\lambda_1+n(\lambda_1+\lambda_2)}{\lambda_1+\lambda_2}
\right).
\label{eq:convexity-ineq}
\end{equation}

\begin{remark}
    Limit set doesn't have to be connected.
\end{remark}

There is the order on the limiting points. Inequality \eqref{eq:convexity-ineq} corresponds to positivity of the rotation between adjacent limiting points; equivalently, the limiting edge is non-negative with respect to the induced lengthwise orientation --- an example of what happens when this inequality is violated is shown on Fig. 21 of \cite{Trivalency}. For strictly convex boundaries, the inequalities \eqref{eq:convexity-ineq} are equalities.

\begin{figure}[H]
\centering
\begin{tikzpicture}[>=Stealth,line cap=round,line join=round,scale=1.28]
  % common stem and first long L/R edges
  \coordinate (S) at (-0.85,-0.68);
  \coordinate (V) at (0,0);
  \coordinate (A) at (0,3.00);
  \coordinate (B) at (2.00,3.00);
  \coordinate (C) at (5.00,0);
  \coordinate (D) at (5.00,0.20);

  \draw[line width=1.15pt] (S)--(V);
  \draw[line width=1.15pt] (V)--(A)--(B);
  \draw[line width=1.15pt] (V)--(C)--(D);

  % one further exact edge on each selected path
  % LRRR... begins with direction (2,-1)
  \coordinate (T) at (3.00,2.50);
  \draw[line width=1.15pt] (B)--(T);

  % RLLLL... begins with direction (-1,2)
  \coordinate (U) at (4.50,1.20);
  \draw[line width=1.15pt] (D)--(U);

  % omitted infinite continuations, shown only by dots
  \coordinate (PL) at (3.50,2.20);
  \coordinate (PR) at (4.20,1.50);
  \foreach \p in {(3.12,2.43),(3.25,2.35),(3.38,2.27)}
    \fill \p circle (0.9pt);
  \foreach \p in {(4.42,1.28),(4.34,1.36),(4.27,1.43)}
    \fill \p circle (0.9pt);
  \fill (PL) circle (1.0pt);
  \fill (PR) circle (1.0pt);

  % a few sibling branches, making the ambient branching visible
  \draw[thin] (A)--++(-0.36,0.36);
  \draw[thin] (B)--++(0,0.42);
  \draw[thin] (C)--++(0.36,-0.36);
  \draw[thin] (D)--++(0.42,0);

  % bold boundary side; with the two selected paths it bounds a convex face
  \draw[line width=2.2pt, gray] (PL)--(PR);
  \node[anchor=west] (lab) at (4.55,2.00) {the boundary side};
  \draw[->,thin] (lab.west)--($(PL)!0.55!(PR)$);
\end{tikzpicture}
\caption{A boundary side is the segment joining the limiting points of two adjacent paths. First: turn left once, then turn always right. Second: turn right once, then turn always left.}
\end{figure}
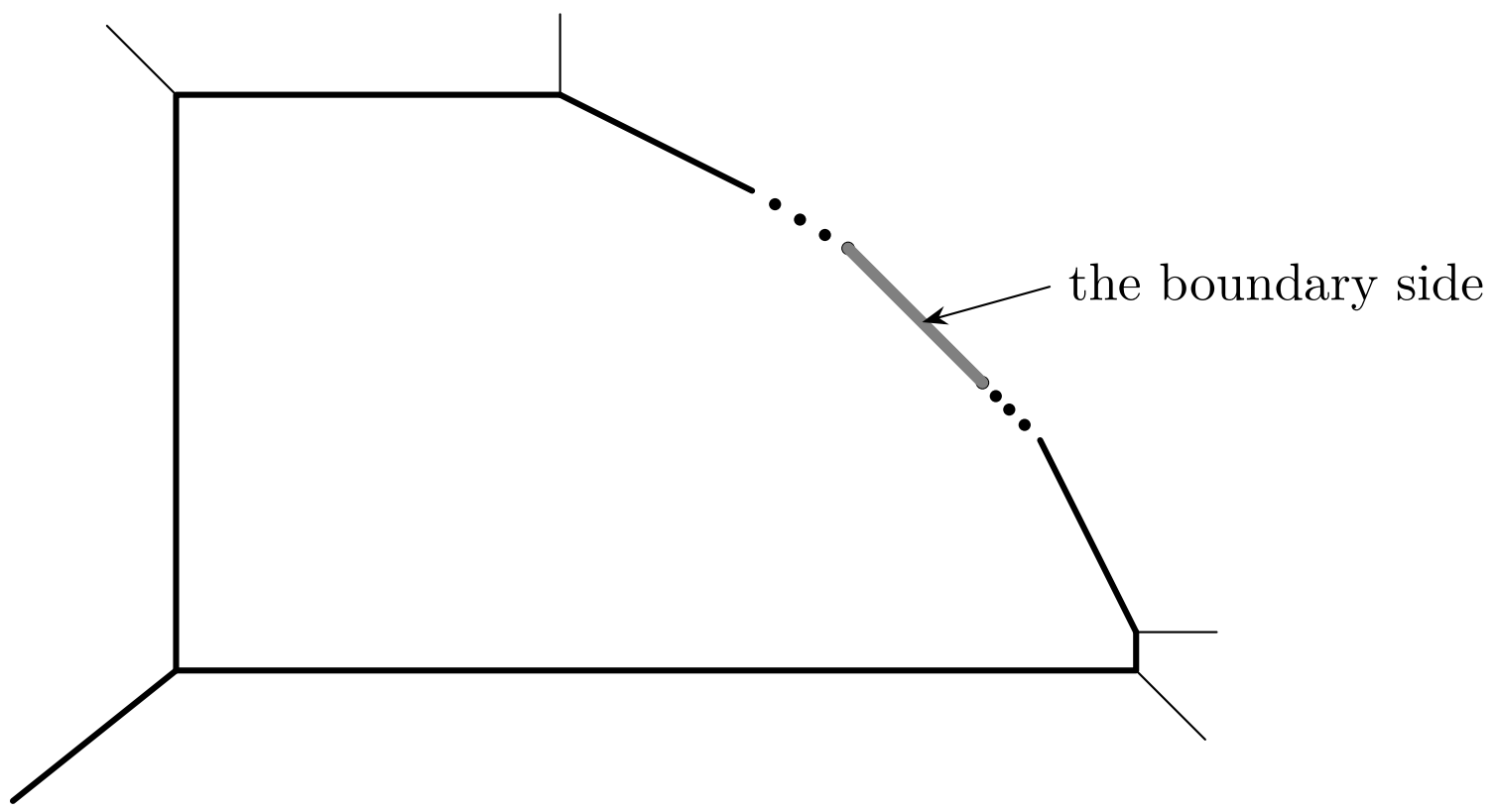

\subsection{Concave limiting curves and the open problem}

Provided that the inequalities~\eqref{eq:convexity-ineq} above hold, for a given $b>0,$ connecting adjacent limiting points, as well as the point $(0,a_{(0,0)})$ with the first limit point, i.e. the result of always going left starting with $(0,0),$  gives a bijective map with plots of non-increasing non-negative concave functions on $[0,b]$ with value $b$ at $0$. Note that in this correspondence $b=a_{(0,0)},$ since all the primitive directions of the edges in the path which always turn left have the second coordinate $1,$ as it is seen from Figure \ref{fig-FareyTriang}, and the sum of potential differences is telescopic.

\textit{Proof sketch (see \cite{KSsummation}).}
For a graph $(x,f(x))$ and a primitive vector
\[
\lambda\in\Z^2_{\ge0},
\]
find, by the Legendre transform, the coefficient $c_\lambda$ of the supporting line
\[
\lambda\cdot z+c_\lambda=0,
\qquad z\in\R^2.
\]%а кто такой z откуда он взялся
Then set $$\potential{\lambda_1}{\lambda_2}=c_{\lambda_1}+c_{\lambda_2}-c_{\lambda_1+\lambda_2}.$$ \\
A branch with such a potential realizes the plot of the given function.

\vspace{0.8em}
\begin{center}
\fbox{\parbox{0.9\textwidth}{
\textbf{Open problem.}
After dropping the inequalities on $\potential{\lambda_1}{\lambda_2}$, describe the set of all limiting curves of generalized tropical caustic branches.\\ Is the map from caustic branches to curves injective?
}}%номера неравенств тут? %inequality~\eqref{eq:convexity-ineq}
\end{center}

\subsection{Remarks}
The above is an elaboration of section 5.3 of \cite{Trivalency}. After the geometry of a single branch is understood, it would be desirable to globalize it, i.e. to understand how the maximal branches can be arranged globally. In the convex domain case,
this is done via symplectic minimal models, with the point type symplectic minimal models being described on Figure 3 of \cite{MIKH}, and the segment type three families given on Figure 7 of \cite{Zeta}. It may happen that convex symplectic minimal models already produce all locally convex closed curves with cusps as possible singularities. In such a case, it would strongly suggest that other approach to overcoming complexity (it is based on lattice approximations with set-theoretic averaging over lattice bases) sketched in 5.2 of \cite{Trivalency} in the case of Jordan curves serving as terminal loci is related to the one above by Legendre duality, and one may hope that corresponding zeta functions, which can be either constructed from the wave front or caustic potential perspective are analytically dependent --- this idea comes from a parallel with gauge zeta functions in contact geometry communicated by Antoine Prouff. The third approach, in this light, presents a mystery, as it deals not with curves but rather with folded domains.

Finally, the higher dimensional case of the above is even more hypothetical. For instance, in dimension two, tropical wave fronts end at points or segments. In dimension three, points and segments also arise as final moments of tropical wave fronts, but we may also see polygons playing such a role. It is not clear which polygons are realized in this way. Also, branching structure is not fully clarified, for instance, we know since very recently that there are critical events in tropical wave fronts which are not blow-downs, but the full classification is missing.

\end{document}